\documentclass[10pt]{amsart}
\usepackage{amsmath,amssymb,amsthm,amscd,mathrsfs}
\usepackage[all]{xy}
\usepackage{color} 
\CompileMatrices
\numberwithin{equation}{section}
\newtheorem{prop}{Proposition}[section]
\newtheorem{theo}[prop]{Theorem}
\newtheorem{lemm}[prop]{Lemma}
\newtheorem{coro}[prop]{Corollary}
\newtheorem{rema}[prop]{Remark}

\newtheorem{defi}[prop]{Definition}

\numberwithin{equation}{section}

\newcommand{\be}{\begin{equation}}
\newcommand{\ee}{\end{equation}}

\newcommand\IZ{\mathbb {Z}}

\newcommand\IQ{\mathbb {Q}}
\newcommand{\IC}{\mathbb{C}}

\newcommand{\IR}{\mathbb{R}}

\newcommand{\ba}{\begin{array}}
\newcommand{\ea}{\end{array}}

\newcommand{\CV}{{\mathcal V}}
\newcommand{\CX}{{\mathcal X}}

\newcommand{\IH}{{\mathbb H}} 
\newcommand{\bal}{\begin{aligned}}
\newcommand{\eal}{\end{aligned}}

\newcommand{\mfd}{{\mathfrak d}}

\newcommand{\mfg}{{\mathfrak g}}

\newcommand{\mfh}{{\mathfrak h}}

\newcommand{\mfm}{{\mathfrak{M}}}

\newcommand{\longto}{\longrightarrow}
\newcommand{\lochom}{{\mathcal Hom}}

\newcommand{\CO}{{\mathcal O}}

\newcommand{\CE}{{\mathcal E}}

\newcommand{\CH}{{\mathcal H}}
\newcommand{\CF}{{\mathcal F}}

\newcommand{\CC}{{\mathcal C}}

\newcommand{\wCC}{{\widetilde {\mathcal C}}}

\newcommand{\mfe}{{\mathfrak e}}

\newcommand{\CQ}{{\mathcal Q}}

\newcommand{\udelta}{{\underline \delta}}

\newcommand{\obj}{{\mathfrak {Ob}}}

\title{Rank Two ADHM Invariants and Wallcrossing}
\author{W.-y. Chuang, D.-E. Diaconescu, G. Pan}

\begin{document}

\begin{abstract} Generalized Donaldson-Thomas invariants corresponding 
to local D6-D2-D0 configurations  
are defined applying the formalism of Joyce and Song to ADHM sheaves on curves. 
A wallcrossing formula for invariants of D6-rank two is proven
and shown to agree with the wallcrossing formula of Kontsevich and 
Soibelman. Using this result, the asymptotic D6-rank two invariants of $(-1,-1)$ and 
$(0,-2)$ local rational curves are computed in terms of the D6-rank one 
invariants. 
\end{abstract} 
\maketitle 

\tableofcontents

\section{Introduction}\label{intro}
Motivated by string theory considerations, 
ADHM invariants of curves were introduced in \cite{modADHM} as an alternative
construction for the local stable pair theory of curves of Pandharipande and Thomas 
\cite{stabpairs-I}.  They have been subsequently generalized in \cite{chamberI} 
employing a natural variation of the 
stability condition. An important feature of this construction resides in its 
compatibility with the Joyce-Song theory of generalized Donaldson-Thomas 
invariants \cite{genDTI}.  Explicit wallcrossing formulas 
for ADHM invariants have been derived and proven in \cite{chamberII}
using Joyce theory \cite{J-I,J-II,J-III,J-IV} and Joyce-Song theory \cite{genDTI}. 

The purpose of the present paper is to study a further generalization 
of ADHM invariants allowing higher rank framing sheaves. 
This generalization is motivated in part 
by recent work of Toda \cite{ranktwo} and Stoppa \cite{ranktwoGW} on rank two 
generalized Donaldson-Thomas invariants of Calabi-Yau threefolds. 
In contrast to \cite{ranktwo}, \cite{ranktwoGW}, the invariants constructed 
here count local objects with nontrivial D2-rank, in physics terminology. 
Similar rank two Donaldson-Thomas invariants of Calabi-Yau threefolds 
are defined and computed in \cite{AS-I, AS-II} using both wallcrossing and 
direct virtual localization methods.  

Local invariants with higher D6-rank are also interesting on physical grounds. 
Explicit results for such invariants are required in order to test the OSV 
conjecture \cite{Ooguri:2004zv} for magnetically charged black holes. 
In particular, such results would be needed in order to extend the work 
of \cite{Aganagic:2004js} to local D-brane configuration with nonzero D6-rank.
According to \cite{DM-split}, counting invariants with higher D6-rank are also 
expected to determine certain subleading 
corrections to the OSV formula \cite{Ooguri:2004zv}. Moreover, 
walls of marginal stability for BPS states  with nontrivial D6-charge in a local 
conifold model have been studied from a supergravity point of view in \cite{JM}. 
The construction presented below should be viewed as a rigorous mathematical 
framework for the microscopic theory of such BPS states. A detailed comparison 
will appear elsewhere. 

From the point of view of six dimensional 
gauge theory dynamics, the invariants constructed in this paper 
can be thought of as a higher rank generalization of local Donaldson-Thomas
invariants of curves. 
It should be noted however that they are not the same as the higher rank 
local DT invariants defined in \cite{modADHM}, which, from a gauge theoretic 
point of view, are Coulomb branch invariants (see also \cite{JM,Cirafici:2008sn} for a noncommutative 
gauge theory approach.) 
Instead, employing a different treatment of boundary conditions in the 
six dimensional gauge theory, the approach presented below yields Higgs branch 
invariants.  

The geometric setup of the present construction is specified by a
 triple $\CX=(X,M_1,M_2)$
where $X$ is a smooth projective curve of $X$ over $\IC$ of genus $g$, and 
$M_1,M_2$ are line bundles on $X$ so that $M=M_1\otimes_X M_2$ 
is isomorphic to the anticanonical bundle $K_X^{-1}$. 
The data $\CX$ determines 
an abelian category $\CC_\CX$ of quiver sheaves on $X$ constructed 
in \cite[Sect 3]{chamberI}. 

Section (\ref{adhmstab}) consists of a step-by-step construction of 
counting invariants for objects of $\CC_\CX$ following \cite{genDTI}. 
The required stability conditions, chamber structure and moduli stacks are 
presented in sections (\ref{defbasic}), (\ref{chamber}), (\ref{moduli}) 
respectively. Some basic homological algebra results are provided in 
section (\ref{extensions}). The construction is concluded in 
section (\ref{invariants}). Given a stability parameter $\delta\in \IR$ 
the geometric data $\CX$
determines a function $A_\delta: \IZ^{\times 3}\to \IQ$, 
which assigns to any triple $\gamma=(r,e,v)$ the virtual number of 
$\delta$-semistable ADHM sheaves on $X$ of type $\gamma$. 
This function is supported on $\IZ_{\geq 1}\times \IZ\times \IZ_{\geq 0}$. 
In physics terms, the integers $(r,e,v)$ correspond to D2, D0 and D6-brane 
charges respectively. In the derivation of wallcrossing formulas, it is more 
convenient to use the alternative notation $\gamma=(\alpha,v)$, 
$\alpha =(r,e) \in \IZ\times\IZ$. Moreover, the invariants $A_\delta(\alpha,0)$ 
are manifestly independent on $\delta$, and will be denoted by $H(\alpha)$
since they are counting invariants for Higgs sheaves on $X$.

Note that for a fixed type $\gamma$ there is a finite set $\Delta(\gamma)\subset 
\IR$ 
of critical stability parameters dividing the real axis in stability chambers
(see lemma (\ref{chamberlemma}). 
The invariants $A_\delta(\gamma)$ are constant when $\delta$ varies within 
a stability chamber. 
The chamber $\delta> \mathrm{max}\, \Delta(\gamma)$ will be referred to as the 
asymptotic chamber, and the corresponding invariants will be also denoted by 
$A_\infty(\gamma)$. 
The main result of this paper is a wallcrossing formula for 
$v=2$ ADHM invariants at a critical stability parameter $\delta_c>0$ 
of type $(\alpha, 2)$, for arbitrary $\alpha =(r,e) 
\in \IZ_{\geq 1}\times \IZ$. Certain preliminary definitions will be needed 
in the formulation of this result, as follows. 

For any integer $l\in \IZ_{\geq 1}$, and any $v\in \{1,2\}$ 
let ${\mathcal {HN}}_-(\alpha,v,\delta_c,l,l-1)$ denote the set of ordered sequences
$((\alpha_i))_{1\leq i\leq l}$, $\alpha_i\in \IZ_{\geq 1}\times \IZ$,
$1\leq i\leq l$ satisfying the following conditions
\be\label{eq:alphadecomp}
\alpha_1+\cdots +\alpha_l=\alpha
\ee
and 
\be\label{eq:slopesA}
{e_1\over r_1} = \cdots = {e_{l-1}\over r_{l-1}}={e_l+v\delta_c\over r_l}
={e+v\delta_c\over r}
\ee
For any integer $l\in \IZ_{\geq 2}$, let ${\mathcal {HN}}_-(\alpha,2,\delta_c,l,l-2)$
denote the set of ordered sequences
$((\alpha_i))_{1\leq i\leq l}$, $\alpha_i\in \IZ_{\geq 1}\times \IZ$,
$1\leq i\leq l$ satisfying condition \eqref{eq:alphadecomp},
\be\label{eq:slopesB}
{e_1\over r_1} = \cdots = {e_{l-2}\over r_{l-2}}={e_{l-1}+\delta_c\over r_{l-1}}
={e_{l}+\delta_c\over r_l} = {e+2\delta_c\over r},
\ee
and 
\be\label{eq:slopesC}
1/{r_{l-1}} < 1/{r_{l}}.
\ee
 Let $0<\delta_-<\delta_c<\delta_+$ be stability parameters so that there are no critical 
 stability parameters of type $(\alpha,2)$ in the intervals $[\delta_-,\ \delta_c)$, 
 $(\delta_c,\ \delta_+]$. For any triple $(\beta,v)$, $\beta\in 
 \IZ_{\geq 1}\times\IZ\times \IZ_{\geq 1}$, $v\in \{1,2\}$, the 
 invariants $A_{\delta_\pm}(\beta,v)$ will be denoted by 
 $A_\pm(\beta,v)$.  
Then the following result holds 
for $\delta_-,\delta_+$ sufficiently close to $\delta_c$.
\begin{theo}\label{wallcrossingthmA}
 The $v=2$ ADHM invariants satisfy the following wallcrossing formula 
\be\label{eq:wallformulaG} 
\bal 
& A_-(\alpha,2)-A_+(\alpha,2) = \\
&\mathop{\sum_{l\geq 2}}_{} {1\over (l-1)!}
\mathop{\sum_{(\alpha_i)\in  {\mathcal{HN}}_-(\alpha,2, \delta_c, l,l-1)}}_{}
A_+(\alpha_l,2)\prod_{i=1}^{l-1}f_2(\alpha_i)H(\alpha_i)\\
&-{1\over 2}\mathop{\sum_{l\geq 1}}_{} {1\over (l-1)!}
\mathop{\sum_{(\alpha_i)\in  {\mathcal{HN}}_-(\alpha,2, \delta_c, l+1,l-1)}}_{}
g(\alpha_{l+1},\alpha_{l}) A_+(\alpha_l,1)
A_+(\alpha_{l+1},1)
\prod_{i=1}^{l-1}f_2(\alpha_i)H(\alpha_i)\\
& +{1\over 2}\mathop{\sum_{(\alpha_1,\alpha_2)\in 
{\mathcal{HN}}_-(\alpha,2,\delta_c,2,0)}}_{} \mathop{\sum_{l_1\geq 1}}_{}
\mathop{\sum_{l_2\geq 1}}_{} {1\over (l_1-1)!}{1\over (l_2-1)!}
\mathop{\sum_{(\alpha_{1,i})\in  {\mathcal{HN}}_-(\alpha_1,1, \delta_c, l_1,l_1-1)}}_{}\\
&
\mathop{\sum_{(\alpha_{2,i})\in  {\mathcal{HN}}_-(\alpha_2,1, \delta_c, l_2,l_2-1)}}_{}
g(\alpha_1,\alpha_2)A_+(\alpha_{1,l_1},1)A_+(\alpha_{2,l_2},1) 
\prod_{i=1}^{l_1-1} f_1(\alpha_{1,i})H(\alpha_{1,i}) 
\prod_{i=1}^{l_2-1} f_1(\alpha_{2,i})H(\alpha_{2,i}) \\
\eal
\ee
where 
\[
\bal
f_v(\alpha) & = (-1)^{v(e-r(g-1))}v(e-r(g-1)),\qquad v=1,2\\
g(\alpha_1,\alpha_2) & = (-1)^{e_1-e_2-(r_1-r_2)(g-1)}(e_1-e_2-(r_1-r_2)(g-1))\\
\eal
\]
for any $\alpha=(e,r)$ respectively $\alpha_i=(r_i,e_i)$, $i=1,2$, 
and the sum in the right hand side of equation \eqref{eq:wallformulaG} 
is finite. 
\end{theo}
Theorem (\ref{wallcrossingthmA}) is proven in section (\ref{wallformula})
using certain stack function identities established in section (\ref{sfidentities}). 
Formula \eqref{eq:wallformulaG} is shown to agree with the wallcrossing 
formula of Kontsevich and Soibelman in section (\ref{KSsect}). 

An application of theorem (\ref{wallcrossingthmA}) to genus zero invariants is  
presented in section (\ref{genuszero}). 
Consider the following generating functions
\be\label{eq:partfctA} 
Z_{\CX,v}(u,q) = \mathop{\sum_{r\geq 1}}_{}\mathop{\sum_{n\in\IZ}}_{} 
u^rq^{n} A_\infty(r,n-r,v)
\ee 
where $v=1,2$. 
Using the wallcrossing formula \eqref{eq:wallformulaG} and the comparison 
result of section (\ref{KSsect}), the following closed formulas
are proven in section (\ref{genuszero}).
\begin{coro}\label{closedform} 
Suppose $X$ is a genus $0$ curve 
and $M_1\simeq \CO_X(d_1)$, $M_2\simeq \CO_X(d_2)$ where 
$(d_1,d_2)=(1,1)$ or $(0,2)$. 
Then 
\be\label{eq:almostproductA}
\bal
&  Z_{\CX,1}(u,q) = \prod_{n=1}^\infty(1-u(-q)^n)^{(-1)^{d_1-1}n}\\
& Z_{\CX,2}(u,q) = {1\over 4}\prod_{n=1}^\infty(1-uq^n)^{2(-1)^{d_1-1}n}
 - {1\over 2}\sum_{\substack{r_1>r_2\geq 1,\ n_1,n_2\in \IZ \\ \mathrm{or}\ 
r_1=r_2\geq1,\ n_2>n_1 \\ \mathrm{or}\ r_1 \geq 1, \ n_1 \in \IZ, \ r_2=n_2=0 }}  (n_1-n_2) (-1)^{(n_1-n_2)} \\
&\qquad \qquad \qquad \qquad A_{\infty}(r_1,n_1-r_1,1)
A_{\infty}(r_2,n_2-r_2,1) u^{r_1+r_2}q^{n_1+n_2}.\\
\end{aligned}
\ee
\end{coro}

\begin{rema} 
The computations in section (\ref{genuszero}) based on the Kontsevich-Soibelman wallcrossing formula can be generalized to invariants of arbitrary rank $v\geq 2$. 
Then it follows that the rank $v$ invariants of local $(-1,-1)$ and 
$(0,-2)$ curves are recursively determined by the 
invariants of lower rank $1\leq v'\leq v$. The resulting formulas are quite complicated, 
and will be omitted. 
\end{rema}

{\it Acknowledgements} 
We are very grateful to Greg Moore for comments and suggestions on the manuscript. 
The work of D.-E. D. is supported in part by NSF grant PHY-0854757-2009. 
WYC is supported by DOE grant DE-FG02-96ER40959.

\section{Higher rank ADHM invariants}\label{adhmstab}

\subsection{Definitions and basic properties}\label{defbasic}
Let $X$ be a smooth projective curve of genus $g\in \IZ_{\geq 0}$ 
over an infinite field $K$ of 
characteristic $0$ equipped with
a very ample line bundle $\CO_X(1)$.  
Let $M_1,M_2$ be fixed line bundles on $X$ equipped with a fixed isomorphism   
$M_1\otimes_X M_2\simeq K_X^{-1}$. Set $M=M_1\otimes_X M_2$.
For fixed data $\CX=(X,M_1,M_2)$,
let $\CQ_{\CX,s}$ denote the abelian category 
of $(M_1,M_2)$-twisted coherent ADHM quiver sheaves. An object 
of $\CQ_\CX$ is given by a collection $\CE=(E,E_{\infty}, 
\Phi_{1}, \Phi_{2}, \phi,\psi)$ 
where 
\begin{itemize} 
\item $E, E_{\infty}$ are coherent $\CO_X$-modules
\item $\Phi_{i}:E\otimes_X M_i \to E$, $i=1,2$ , $
\phi:E\otimes_X M_1\otimes_X M_2 \to E_{\infty}$, $ 
\psi:E_{\infty}\to E$ are morphisms of $\CO_X$-modules 
satisfying the ADHM relation 
\be\label{eq:ADHMrelation}
\Phi_1\circ(\Phi_2\otimes 1_{M_1}) - \Phi_2\circ (\Phi_1\otimes 1_{M_2}) 
+ \psi\circ \phi =0.
\ee
\end{itemize}
The morphisms are natural morphisms of quiver sheaves i.e. collections 
$(\xi,\xi_{\infty}):(E,E_{\infty}) \to 
(E',E'_{\infty})$ 
of morphisms of $\CO_X$-modules satisfying the obvious compatibility conditions 
with the ADHM data. 

Let $\CC_{\CX}$ be the full abelian 
subcategory of $\CQ_{\CX}$ consisting of objects with 
$E_{\infty} = V \otimes \CO_X$, where $V$
is a  finite dimensional  
vector spaces over $K$ (possibly trivial.) 
Note that given any two objects $\CE,\CE'$ 
of $\CC_\CX$, the morphisms 
$\xi_{\infty}: V\otimes \CO_X\to V'\otimes \CO_X$ must be of the form 
$\xi_{\infty} = f\otimes 1_{\CO_X}$, where $f:V\to V'$
is a  linear map. 

An  object $\CE$ of $\CC_{\CX}$ will be called locally free if $E$ is a coherent 
locally free $\CO_X$-module. 
Given a coherent $\CO_X$-module $E$  we will denote by $r(E)$, $d(E)$, 
$\mu(E)$ the rank, degree, respectively slope of $E$ if $r(E)\neq 0$. 
The type of an object $\CE$ of $\CC_{\CX}$ is the collection 
$(r(\CE), d(\CE), v(\CE))= (r(E),d(E),\mathrm{dim}(V)))\in 
\IZ_{\geq 0}\times \IZ\times \IZ_{\geq 0}$.
An object of $\CO_X$ will be called an ADHM sheaf in the following. 
Throughout this paper, the integer $v(\CE)$ will be called the rank of $\CE$, 
as opposed to the terminology used in \cite{modADHM,chamberI,chamberII}, 
where the rank of $\CE$ was defined to be $r(\CE)$. 
Note that the objects of $\CC_{\CX}$ with $v(\CE)=0$ 
form a full abelian category which is  naturally equivalent
to the 
abelian category of
 Higgs sheaves 
on $X$ with coefficient bundles $(M_1,M_2)$ (see for example 
\cite[App. A]{chamberI} for brief summary of the relevant definitions.)

Let 
$\delta\in \IR$ be a stability parameter. 
The $\delta$-degree of an object $\CE$ of $\CC_\CX$ is defined by 
\be\label{eq:udeltadeg} 
\mathrm{deg}_\udelta(\CE) = d(\CE) +  \delta v(\CE).
\ee
If $r(\CE)\neq 0$, the $\delta$-slope of $\CE$ is defined by 
\be\label{eq:udeltaslope}
\mu_\delta(\CE) = {\mathrm{deg}_{\delta}(\CE)\over r(\CE)}.
\ee
\begin{defi}\label{udeltastability}
Let $\delta\in \IR$ be a stability parameter. 
A nontrivial object $\CE$ of $\CC_\CX$ is  
$\delta$-(semi)stable if 
\be\label{eq:udeltastabcondA} 
r(E)\, \mathrm{deg}_{\delta}(\CE') \ (\leq) \ r(E')\, \mathrm{deg}_\delta(\CE)
\ee
for any proper nontrivial subobject $0\subset \CE'\subset \CE$. 
\end{defi} 

The following lemmas summarize some basic properties of $\delta$-semistable ADHM sheaves.  
The proofs are either standard or very similar to those of \cite[Lemm. 2.4]{chamberI}, \cite[Lemm 3.7]{chamberI} and will be omitted. 
\begin{lemm}\label{basicprop} 
Suppose $\CE$ is a $\delta$-semistable framed ADHM sheaf with $r(\CE)>0$
for some $\delta \in \IR$. Then 
\begin{itemize}
\item[$(i)$] $E$ is locally free. 
\item[$(ii)$]   If $\delta>0$, there is no nontrivial 
linear subspace $0\subset V' \subseteq V$ so that $\psi|_{V'\otimes 
\CO_X}$ is identically zero. Similarly, if $\delta<0$, 
 there is no proper linear subspace 
$0\subseteq V'\subset V$ so that $\mathrm{Im}
(\phi)\subseteq V'\otimes \CO_X$.
  \item[$(iii)$] 
 If $\CE$ is $\delta$-stable
any endomorphism of $\CE$ in $\CC_{\CX}$ is either 
trivial or an isomorphism. If the ground field $K$ is algebraically closed, the endomorphism 
ring of $\CE$ is canonically isomorphic to $K$. 
\end{itemize} 
\end{lemm}

\begin{lemm}\label{boundlemma}
For fixed $(r,e,v)\in \IZ_{>0}\times \IZ\times \IZ_{\geq 0}$ 
there is a constant $c\in \IR$ (depending only on $\CX$ and (r,e,v)) 
so that for any $\delta\in \IR$, 
any $\delta$-semistable framed 
ADHM sheaf of type $(r,e,v)$ satisfies 
\[
\mathrm{\mu_{\max}(E)} < c. 
\]
In particular, the set of isomorphism classes of framed ADHM sheaves of fixed type 
$(r,e,v)$ 
which are $\delta$-semistable for some $\delta\in \IR$
is bounded. 
\end{lemm} 

Given a locally free ADHM sheaf  $\CE=(E, \Phi_1,\Phi_2,\phi,\psi)$ on $X$
of type $(r,e,v)\in \IZ_{\geq 1}\times \IZ\times \IZ_{\geq 0}$,
the data 
\be\label{eq:dualADHM} 
\begin{aligned}
{\widetilde E} & = E^\vee \otimes_X M^{-1}\\
{\widetilde \Phi}_i & = (\Phi_i^\vee\otimes 1_{M_i}) \otimes 1_{M^{-1}}: {\widetilde E}
\otimes 
M_i \to {\widetilde E} \\
{\widetilde \phi} & = \psi^\vee \otimes 1_{M^{-1}} : {\widetilde E}\otimes_X {M} \to 
V^\vee\otimes \CO_X
\\
{\widetilde \psi} & = \phi^\vee : V^\vee \otimes \CO_X \to {\widetilde E} \\
\end{aligned}
\ee
with $i=1,2$, determines a locally free ADHM sheaf ${\widetilde \CE}$ of type 
$(r,-e+2r(g-1),v)$ where $g$ is the 
genus of $X$. ${\widetilde \CE}$ will be called the dual of $\CE$ in the following. 
Then the following lemma is straightforward.

\begin{lemm}\label{duallemma}
Let $\delta\in \IR$ be a  stability parameter and let 
$\CE$ be a locally free ADHM sheaf on $X$. 
Then $\CE$ is $\delta$-(semi)stable if and only if ${\widetilde \CE}$ is 
$(-\delta)$-(semi)stable.
\end{lemm}

\subsection{Chamber structure}\label{chamber}
This subsection summarizes the main properties of 
$\delta$-stability chambers.
\begin{defi}\label{asympstab} 
An ADHM sheaf $\CE$ of type 
$(r,e,v)\in\IZ_{\geq 1}\times \IZ\times \IZ_{\geq 0}$ is 
asymptotically (semi)stable if the following 
conditions hold 
\begin{itemize} 
\item[$(i)$] $E$ is locally free, $\psi:V\otimes \CO_X\to E$ is not 
identically zero, and 
there is no saturated proper nontrivial subobject $0\subset \CE'\subset \CE$ 
in $\CC_\CX$ so that $v(\CE')/r(\CE') > v/r$. 
\item[$(ii)$] Any proper nontrivial subobject $0\subset \CE'\subset \CE$ 
with $v(\CE')/r(\CE') = v/r$ satisfies the slope inequality $\mu(E') \ (\leq) \ \mu(E)$. 
\end{itemize} 
\end{defi} 

Here a subobject $\CE'\subset \CE$ is called saturated in the underlying 
coherent sheaf $E'$ is saturated in $E$. Note that according to \cite[Lemm. 
3.10]{chamberI}, any proper subobject $0\subset \CE'\subset \CE$ admits a canonical saturation ${\overline {\CE'}}\subset \CE$.

\begin{lemm}\label{asympbound}
The set of isomorphism classes of asymptotically semistable ADHM sheaves 
of fixed type $(r,e,v)\in\IZ_{\geq 1}\times \IZ\times \IZ_{\geq 1}$ is bounded. 
\end{lemm}

{\it Proof.} The proof is based on Maruyama's boundedness theorem.
Suppose $\CE$ is asymptotically semistable of type $(r,e,v)$, and the underlying 
coherent sheaf $E$ is not semistable. Then there is a nontrivial Harder-Narasimhan 
filtration 
\[
0\subset E_1\subset \cdots \subset E_h = E
\]
with $h\geq 2$ so that $\mu(E_j) > \mu(E)$ and $r(E_j)<r$ for all $1\leq j \leq h-1$. 
Suppose $E_j$ is $\Phi_i$-invariant, $i=1,2$, and $\mathrm{Im}(\psi)\subseteq E_j$ 
for some $1\leq j\leq h-1$. Then the data $\CE_j=(E_j, \Phi_i|_{E_j\otimes_X M_i}, 
\phi|_{E_j\otimes_X M}, \psi)$ is subobject of 
$\CE$ with 
\[
v(\CE_j)/r(\CE_j) = {v\over r(\CE_j)}>{v\over r}.
\]
Since $E_j\subset E$ is saturated, it follows that $\CE_j$ violates condition $(i)$
in definition (\ref{asympstab}). Therefore for any $1\leq j\leq h$, $E_j$ is either 
not preserved by some $\Phi_i$, $i=1,2$, or it does not contain the image 
of $\psi$. From this point on the proof is identical to the proof of 
\cite[Prop. 2.7]{modADHM}. 

\hfill $\Box$

\begin{defi}\label{genstab} 
Let $\delta\in \IR_{>0}$. 
A $\delta$-semistable ADHM sheaf $\CE$ of type 
$(r,e,v)\in \IZ_{\geq 1}\times \IZ\times \IZ_{\geq 0}$ is generic 
if it is either $\delta$-stable or any proper nontrivial subobject $0\subset \CE'\subset \CE$ 
of type $(r',e',v')\in \IZ_{\geq 1}\times \IZ\times \IZ_{\geq 0}$ satisfies 
\be\label{eq:gencond} 
{e'\over r'}={e\over r} \qquad {v'\over r'} = {v \over r}.
\ee
The stability parameter $\delta\in \IR_{>0}$ is called generic of type $(r,e,v)$ if 
any $\delta$-semistable ADHM sheaf of type $(r,e,v)$ is generic. The stability 
parameter $\delta \in \IR_{>0}$ is called critical of type $(r,e,v)$ 
if there exists a nongeneric 
$\delta$-semistable ADHM sheaf of type $(r,e,v)$. 
\end{defi}

Lemma (\ref{boundlemma}) implies the following.
\begin{lemm}\label{asympstablemma} 
For fixed $(r,e,v)\in \IZ_{\geq 1}\times \IZ\times \IZ_{\geq 1}$ 
there exists $\delta_\infty\in\IR_{>0}$ so that  
for all $\delta\geq \delta_\infty$  
an ADHM sheaf $\CE$ of type $(r,e,v)$ is $\delta$-(semi)stable 
if and only if it is asymptotically (semi)stable. 
\end{lemm}

{\it Proof.} The proof if similar to the proof of lemma \cite[Lemm. 4.7]{chamberI}. 
Some details will be provided for convenience. It is straightforward to prove that 
asymptotic stability implies $\delta$-stability for sufficiently large $\delta$ 
using lemma (\ref{boundlemma}). The converse is slightly more involved. 
First note that given any nontrivial locally free ADHM sheaf $\CE$, any linear subspace 
$V'\subset V$, determines a canonical subobject $\CE_{V'}\subset \CE$.
$\CE_{V'}$ is the saturation of the subobject of $\CE$ generated by 
$V'\otimes \CO_X$ by successive applications of the ADHM morphisms $\psi, \Phi_i, \phi$. 
Since $\CE_{V'}$ is canonically determined by $V'$ and $\CE$, lemma 
(\ref{boundlemma}) implies that the set of isomorphism classes of subobjects 
$\CE_{V'}$, where $\CE$ is a $\delta$-semistable ADHM sheaf of type 
$(r,e,v)$ for some $\delta>0$ is bounded. Moreover, by construction, any 
subobject $0\subset \CE'\subset \CE$ contains the canonical subobject $\CE_{V'}$. 

Now suppose that for any $\delta>0$ there exists a $\delta$-semistable ADHM 
sheaf $\CE$ of type $(r,e,v)$ which is not asymptotically stable. 
Let $0\subset \CE'\subset \CE$ be a saturated nontrivial proper saturated subobject 
violating the asymptotic stability conditions. Note that $\CE'$ cannot violate 
condition $(ii)$ in definition (\ref{asympstab}) since $\CE$ is $\delta$-semistable. 
Therefore it must violate condition $(i)$ i.e. $v'/r'>v/r$ where 
$r'=r(\CE')$. In particular $v'=v(\CE')>0$. 
Then the subobject $\CE_{V'}$ also violates condition $(i)$ since 
\[
{v(\CE_{V'})\over r(\CE_{V'})} = {v'\over r(\CE_{V'})} \geq {v'\over r'}>{v/r}.
\]
Since $\CE$ is $\delta$-semistable $\mu_\delta(\CE_{V'}) \leq 
\mu_\delta(\CE)$. However, as noted above, the set of isomorphism classes 
of all $\CE_{V'}$ is bounded, therefore the set of all types $(r(\CE_{V'}), d(\CE_{V'}), 
v(\CE_{V'}))$ is finite. Taking $\delta$ sufficiently large, this leads to a contradiction. 

\hfill $\Box$

By analogy with \cite[Lemm. 4.4]{chamberI}, 
\cite[Lemm. 4.6]{chamberI}, lemmas (\ref{asympstablemma}) and 
(\ref{duallemma}) imply the following. 

\begin{lemm}\label{chamberlemma}
Let $(r,e,v)\in\IZ_{\geq 1}\times \IZ\times \IZ_{\geq 1}$ be a fixed type. 
Then there is a finite set $\Delta(r,e,v) \subset \IR$ of critical stability parameters of type 
$(r,e,v)$. 
Given any two stability parameters $\delta,\delta'\in \IR$, $\delta< \delta'$
so that $[\delta, \ \delta']\cap \Delta(r,e,v)=\emptyset$, the set of $\delta$-semistable 
ADHM sheaves of type $(r,e,v)$ is identical to the set of $\delta'$-semistable 
ADHM sheaves of type $(r,e,v)$.
\end{lemm}

\begin{rema}\label{requalsone}
It is straightforward to check that $\Delta(1,e,v)=\{0\}$ for any $v\geq 1$. 
\end{rema} 

\begin{lemm}\label{specstab} 
Let $(r,e,v)\in \IZ_{\geq 1} \times \IZ\times \IZ_{\geq 1}$ 
and let $\delta_c>0$ be a critical stability parameter of type $(r,e,v)$.
Let $\delta_\pm >0$ be stability parameters so that $\delta_-<\delta_c<\delta_+$ 
and $[\delta_-,\ \delta_c)\cap \Delta(r,e,v) =\emptyset$, 
$(\delta_c,\ \delta_+]\cap \Delta(r,e,v) =\emptyset$. 
If $\CE$ is a $\delta_\pm$-semistable ADHM sheaf of type $(r,e,v)$, then 
$\CE$ is also $\delta_c$-semistable. 
\end{lemm}

\begin{defi}\label{admconfig}
Let $(r,v)\in \IZ_{\geq 1} \times \IZ_{\geq 1}$. 

$(a)$ A positive admissible configuration of type $(r,v)$
is an ordered sequence of integral points 
$\left(\rho_i=(r_i,v_i)\in \IZ_{\geq 1} \times \IZ_{\geq 0}\right)_{1\leq i\leq h, \, h\geq 1}$ satisfying the following conditions 
\begin{itemize}
\item $\rho_1+\cdots + \rho_h = (r,v)$.
\item $(v_1+\cdots +v_i)/(r_1+\cdots+r_i) > v/r$ 
and $v_{i}/r_i > v_{i+1}/r_{i+1}$ for all $i=1,\ldots, h-1$.
\end{itemize}
$(b)$ A negative admissible configuration of type $(r,v)$
is an ordered sequence of integral points 
$\left(\rho_i=(r_i,v_i)\in \IZ_{\geq 1} \times \IZ_{\geq 0}\right)_{1\leq i\leq h,\, h\geq 1}$ satisfying the following conditions 
\begin{itemize}
\item 
$\rho_1+\cdots + \rho_h = (r,v)$.  
\item $(v_1+\cdots +v_i)/(r_1+\cdots+r_i) < v/r$ 
and $v_{i}/r_i< v_{i+1}/r_{i+1}$ for all $i=1,\ldots, h-1$.
\end{itemize}\end{defi}

\begin{rema}\label{configrema}
$(i)$ 
It is straightforward to prove that for fixed $(r,v)\in \IZ_{\geq 1} \times \IZ_{\geq 1}$ the set of positive, respectively negative, admissible configurations 
is finite. These sets will be denoted by ${\mathcal {HN}}_\pm(r,v)$.

$(ii)$ The only positive, respectively negative admissible configuration of type $(r,v)$ 
with $h=1$ is $(\rho=(r,v))$.
\end{rema}

\begin{lemm}\label{HNfiltrations} 
Let $\delta_c\in \IR_{>0}$ be a critical stability parameter of 
type $(r,e,v)\in \IZ_{\geq 1}\times \IZ\times \IZ_{\geq 1}$. 
Then the following hold. 

$(i)$ There exists $\epsilon_+ >0$, 
so that $(\delta_c,\ \delta_c+\epsilon_+]\cap \Delta(r,e,v)=\emptyset$ 
and the following holds for any $\delta_+\in (\delta_c, \ \delta_c+\epsilon_+)$. 
A locally free ADHM sheaf $\CE$ of type $(r,e,v)$ 
on $X$ is $\delta_c$-semistable if and only if 
it is either $\delta_+$-semistable or there exists a unique filtration of the form 
\be\label{eq:plusfiltr} 
0= \CE_0\subset \CE_1\subset \cdots \subset \CE_h=\CE 
\ee
with $h\geq 2$ satisfying the following conditions 
\begin{itemize}
\item The successive quotients $\CF_i=\CE_i/\CE_{i-1}$, $i=1,\ldots, h$ of the 
filtration \eqref{eq:plusfiltr} are locally free ADHM sheaves with numerical types 
$(r_i,e_i,v_i)\in \IZ_{\geq 1}\times \IZ\times \IZ_{\geq 0}$. 
$\delta_+$ is noncritical of type $(r_i,e_i,v_i)$, $\CF_i$ is 
$\delta_+$-semistable and $\mu_{\delta_c}(\CF_i) = \mu_{\delta_c}(\CE)$
for all $i=1,\ldots, h$.
\item The sequence $\rho_i=(r_i,v_i)$, $i=1,\ldots, h$ 
is a positive admissible configuration of type $(r,e,v)$. 
\end{itemize} 

$(ii)$ There exists $\epsilon_- >0$, 
so that $[\delta_c-\epsilon_-, \delta_c)\cap \Delta(r,e,v)=\emptyset$ 
and the following holds for any $\delta_-\in (\delta_c-\epsilon_-, \ \delta_c)$. 
A locally free ADHM sheaf $\CE$ of type $(r,e,v)$ 
on $X$ is $\delta_c$-semistable if and only if 
it is either $\delta_-$-semistable or there exists a unique filtration of the form 
\be\label{eq:minusfiltr} 
0= \CE_0\subset \CE_1\subset \cdots \subset \CE_h=\CE 
\ee
with $h\geq 2$ satisfying the following conditions 
\begin{itemize}
\item The successive quotients $\CF_i=\CE_i/\CE_{i-1}$, $i=1,\ldots, h$ of the 
filtration \eqref{eq:minusfiltr} are locally free ADHM sheaves with numerical types 
$(r_i,e_i,v_i)\in \IZ_{\geq 1}\times \IZ\times \IZ_{\geq 0}$. 
$\delta_-$ is noncritical of type $(r_i,e_i,v_i)$, $\CF_i$ is 
$\delta_-$-semistable and $\mu_{\delta_c}(\CF_i) = \mu_{\delta_c}(\CE)$
for all $i=1,\ldots, h$.
\item The sequence $\rho_i=(r_i,v_i)$, $i=1,\ldots, h$ 
is a negative admissible configuration of type $(r,e,v)$. 
\end{itemize} 
\end{lemm}

{\it Proof.} 
The proof is similar to the proof of \cite[Lemm. 4.13]{chamberI}. 
Details are included below for completeness. 
Note that it suffices to prove statement $(i)$ since the proof of $(ii)$ is 
analogous. 

Let $\delta_+>\delta_c$ be an arbitrary 
 noncritical stability parameter 
of type $(r,e,v)$ so that $(\delta_c,\delta_+]\cap \Delta(r,e,v)=\emptyset$.   
Suppose $\CE$ is a $\delta_c$-semistable ADHM 
sheaf on $X$. Then $\CE$ is either 
$\delta_+$-stable or there is a Harder-Narasimhan filtration of $\CE$ with respect to 
$\delta_+$-semistability
\be\label{eq:plusHNfiltr}
0\subset \CE_1 \subset \cdots \subset \CE_h = \CE 
\ee
where $h\geq 2$. It is straightforward to check that $\CE_l$, $1\leq l\leq h$ 
must have $r(\CE_l)\geq 1$ and the successive quotients 
$\CF_l$, $0\leq l \leq h-1$ must also have $r_l\geq 1$.  
Then by the general properties of Harder-Narasimhan filtrations 
\be\label{eq:HNineqB}
\mu_{\delta_+}(\CE_1)> \mu_{\delta_+}(\CE_2/\CE_1)>\cdots >
\mu_{\delta_+}(\CE_h/\CE_{h-1})
\ee
and 
\be\label{eq:ineqAB}
\mu_{\delta_+}(\CE_l)> \mu_{\delta_+}(\CE)
\ee
for all $1\leq l\leq h-1$. Since $\CE$ is $\delta_c$-semistable by assumption, 
inequalities \eqref{eq:ineqAB} imply that 
\be\label{eq:HNineqX}
v(\CE_l)/r(\CE_l) > v/r
\ee for all $l=1,\ldots, h$. Note that $v(\CE_l) = v_1+\cdots+ v_l$, 
$r(\CE_l)=r_1+\cdots + r_l$ for any $l=1,\ldots, h$. 

Moreover, using the $\delta_c$-semistability condition and 
inequalities \eqref{eq:ineqAB} 
we have 
\be\label{eq:ineqB}
\delta_+\left({v\over r}-{v(\CE_l)\over r(E_l)}\right) 
< \mu(E_l) -\mu(E) \leq \delta_c\left({v\over r}-{v(\CE_l)\over r(E_l)}\right)
\ee
for all $l=1,\ldots, h$. 

Now let $\gamma >\delta_c$ be a fixed stability 
parameter so that  $(\delta_c,\gamma]\cap \Delta(r,e,v)=\emptyset$.   
Using Grothendieck's lemma and lemma (\ref{boundlemma}), inequalities \eqref{eq:ineqB} 
 imply that the set of isomorphism classes of locally free ADHM sheaves $\CE'$ on $X$ satisfying condition $(\star)$ below is bounded.
\begin{itemize}
\item[$(\star)$] 
There exists a $\delta_c$-semistable 
ADHM sheaf $\CE$ of type $(r,e,v)$ and a stability parameter $\delta_+
\in (\delta_c\ \gamma]$ so that $\CE'\simeq \CE_l$ for 
some $l\in \{1,\ldots, h\}$, where $0\subset \CE_1\subset\cdots\subset 
\CE_h=\CE$, 
$h\geq 1$, is the  Harder-Narasimhan filtration of $\CE$ 
with respect to $\delta_+$-semistability. 
\end{itemize} 
Then it follows that the set 
of numerical types $(r',e',v')$ of 
locally free ADHM sheaves $\CE'$ satisfying  property $(\star)$ 
is finite. This implies that there exists $0< \epsilon_+<\gamma-\delta_c$ 
so that for any $\delta_+\in (\delta_c, \ \delta_c+\epsilon_+)$, and 
any $\delta_c$-semistable ADHM sheaf $\CE$ of type $(r,e,v)$
inequalities \eqref{eq:ineqB} can be satisfied only if 
\be\label{eq:equal}
\mu_{\delta_c}(\CE_l) = \mu_{\delta_c}(\CE)
\ee
for all $l=1,\ldots, h$. 
Hence also 
\[ 
\mu_{\delta_c}(\CE_l/\CE_{l-1}) = \mu_{\delta_c}(\CE)
\]
for all $l=2,\ldots, h$. Then inequalities \eqref{eq:HNineqB}, \eqref{eq:HNineqX}
 imply 
that the sequence $\rho_l=(r_l,v_l)$, $l=1,\ldots, h$ is a 
positive admissible configuration. 
Therefore for all $\delta_+\in (\delta_c, \ \delta_c+\epsilon_+)$, 
any locally free $\delta_c$-semistable ADHM sheaf $\CE$ of type $(r,e,v)$ is either 
$\delta_+$-stable or has a 
Harder-Narasimhan filtration with respect to $\delta_+$-semistability 
as in lemma (\ref{HNfiltrations}.$i$). 

Next note that the set of numerical types 
\be\label{eq:slopeset} 
{\sf S}_{\delta_c}(r,e,v) = \{(r',e',v')\in \IZ_{\geq 1}\times \IZ\times \IZ_{\geq 0}\, |\, 0<r'\leq r ,\ 0\leq v'\leq v, \
r(e'+\delta_c v') = r'(e+\delta_c v)\}
\ee
is finite. Therefore $0< \epsilon_+<\gamma-\delta_i$ above may be chosen so 
that there are no critical stability parameters of type $(r',e',v')$ in the interval 
$(\delta_c, \ \delta_c+\epsilon_+)$ for any
$(r',e',v')\in {\sf {S}}_{\delta_c}(r,e,v)$. In particular, $\delta_+$ is noncritical of type 
$(r_i,e_i,v_i)$, $i=1,\ldots, h$ for any Harder-Narasihan filtration as above. 

Conversely, suppose $\CE$ is a locally free ADHM sheaf of type $(r,e,v)$ 
on $X$ which has a 
filtration of the form \eqref{eq:plusfiltr} with $\CE'$ $\delta_+$-stable 
and satisfying the conditions of lemma (\ref{HNfiltrations}.$i$) for some 
$\delta_+\in (\delta_c, \ \delta_c+\epsilon_+)$. 
By the above choice of $\epsilon_+$, there are no critical stability parameters 
of type $(r_i,e_i,v_i)$ in the interval $(\delta_c, \ \delta_c+\epsilon_+)$, 
for any $i=1,\ldots,h$. Since $\CF_i$ are $\delta_+$-semistable, lemma (\ref{specstab}) 
implies that $\CF_i$ is also $\delta_c$-semistable, for any $i=1,\ldots,h$. 
Hence $\CE$ is also $\delta_c$-semistable since the $\CF_i$ have equal $\delta_c$-slopes. 

\hfill $\Box$

\subsection{Extension groups}\label{extensions}
Let $\CE',\CE''$ be nontrivial locally free objects in $\CC_\CX$ of types 
$(r',e',v'), (r'',e'',v'')\in \IZ_{\geq 1}\times \IZ\times \IZ_{\geq 0}$. 
Let $\CC(\CE'',\CE')$ be the three term complex 
\be\label{eq:hypercohA} 
\begin{aligned} 
0 \to \begin{array}{c} \lochom_{X}(E'',E') \\ \end{array}
& {\buildrel d_1\over \longto} 
\begin{array}{c}  \lochom_{X}(E''\otimes_{X}M_1 ,E') \\ \oplus \\
\lochom_{X}(E''\otimes_{X} M_2, E') \\ \oplus \\ 
\lochom_{X}(E''\otimes_{X} M,V'\otimes \CO_X)\\ \oplus \\ 
\lochom_{X}(V''\otimes \CO_X,E') \\ \end{array} 
 {\buildrel d_2\over \longto} 
\lochom_{X}(E''\otimes_{X}M,E') \to 0 \\
\end{aligned}
\ee
where 
\[
\begin{aligned}
d_1(\alpha) = 
(& -\alpha \circ \Phi_{1}'' +\Phi_{1}'\circ (\alpha\otimes 1_{M_1}), 
-\alpha \circ \Phi_{2}''+\Phi_{2}'\circ (\alpha\otimes 1_{M_2}),\\ 
&  \phi' \circ (\alpha
\otimes 1_M), 
-\alpha\circ \psi'' )\\
\end{aligned}  
\]
for any local sections $(\alpha,\alpha_\infty)$ of 
the first term 
and 
\[
\begin{aligned}
d_2(\beta_1,\beta_2,\gamma, \delta) = &
\beta_1 \circ (\Phi''_2\otimes 1_{M_1}) -
\Phi_{2}'\circ (\beta_1\otimes 1_{M_2}) 
- \beta_2\circ (\Phi''_{1}\otimes 1_{M_2})\\
 & + \Phi_{1}'\circ (\beta_2\otimes 1_{M_1}) + 
\psi'\circ \gamma + \delta \circ \phi''\\
\end{aligned}
\]
for any local sections $(\beta_1,\beta_2,\gamma, \delta)$ 
of the middle term. The degrees of the three terms in 
\eqref{eq:hypercohA} are $0,1,2$ respectively. 

Let $C(\CC(\CE'',\CE'))$ be the double complex obtained from $\CC(\CE'',\CE')$
by taking 
${\check{\rm C}}$ech resolutions and let ${D}(\CE',\CE'')$ be the diagonal complex 
of $C(\CC(\CE'',\CE'))$. Note that there is a canonical linear map 
\[
\bal 
\mathrm{Hom}(V'',V') & \to D^1(\CE',\CE'') = C^0(\CC^1(\CE'',\CE'))\oplus 
C^1(\CC^0(\CE'',\CE'))\\
f & \to \left[\begin{array}{c} 
{}^t(0,0,-(f\otimes 1_{\CO_X})\circ \phi'', \psi'\circ (f\otimes 1_{\CO_X}))\\
0\\ \end{array}\right]\\
\eal
\]
Given the above expressions for the differentials $d_1,d_2$ it is straightforward 
to check that this map yields a morphism of complexes 
\[
\bal
\varrho: \mathrm{Hom}(V'',V')[-1] & \to {D}(\CE'', \CE') \\
\eal
\]
Let ${\widetilde D}(\CE'',\CE')$ denote the cone of $\varrho$. 
Then the lemma below follows either by explicit $\check{\rm{C}}$ech cochain computations 
as in \cite[Sect. 4]{modADHM} or using the methods of \cite{homquiv}. 
\begin{lemm}\label{extlemma} 
The extension groups $\mathrm{Ext}^k_{\CC_\CX}(\CE'',\CE')$, $k=0,1$ 
are isomorphic to the cohomology groups $H^k({\widetilde D}(\CE'',\CE'))$, $k=0,1$. 
Moreover there is an exact sequence 
\be\label{eq:extseq}
\xymatrix{
0\ar[r] & \IH^0(\CC(\CE'',\CE')) \ar[r] & \mathrm{Ext}^0_{\CC_\CX}(\CE'',\CE')\ar[r]
& {Hom}(V'',V') \\
  \ar[r] & \mathrm{Ext}^1_{\CC_\CX}(\CE'',\CE') \ar[r] &
 \IH^1(\CC(\CE'',\CE'))  \ar[r] & 0  \\}
 \ee
 where $\IH^k(\CC(\CE'',\CE'))$, $k=0,1$ are hypercohomology groups of the 
 complex $\CC(\CE'',\CE')$.
\end{lemm}

\begin{coro}\label{extcoro}
Given any two locally free objects $\CE',\CE''$ 
\be\label{eq:eulerchar} 
\bal
& \mathrm{dim}(\mathrm{Ext}^0_{\CC_\CX}(\CE'',\CE')) -
\mathrm{dim}(\mathrm{Ext}^1_{\CC_\CX}(\CE'',\CE')) - 
\mathrm{dim}(\mathrm{Ext}^0_{\CC_\CX}(\CE',\CE'')) \\  & +
\mathrm{dim}(\mathrm{Ext}^1_{\CC_\CX}(\CE',\CE'')) =
v'e''-v''e'-(v'r''-v''r')(g-1) \\
\eal
\ee
\end{coro}

{\it Proof.} Follows from the exact sequence \eqref{eq:extseq} and 
the fact that the hypercohomology groups of the complex 
$\CC(\CE'',\CE')$ satisfy the duality relation 
\[
\IH^k(\CC(\CE'',\CE')) \simeq \IH^{3-k}(\CC(\CE',\CE''))^\vee
\]
for $k=0,\ldots,3$.

\hfill $\Box$

\subsection{Moduli stacks}\label{moduli}
In the following let the ground field $K$ be $\IC$. 
Let ${\mathfrak {Ob}}(\CX)$ denote the moduli stack of all objects of the abelian 
category  $\CC_\CX$ and let ${\mathfrak {Ob}}(\CX,r,e,v)$ denote the open 
and closed component of type $(r,e,v)\in \IZ_{\geq 1}\times \IZ\times \IZ_{\geq 0}$. 
Standard arguments analogous to \cite[Sect. 9]{J-I}, \cite[Sect. 10]{J-I}
prove that  ${\mathfrak {Ob}}(\CX)$ is an algebraic 
stack locally of finite type and it satisfies conditions \cite[Assumption 7.1]{J-I}, 
\cite[Assumption 8.1]{J-I}. Given the boundedness result 
(\ref{boundlemma}), the  following is also standard. 
\begin{prop}\label{modstack} 
For fixed type $(r,e,v)\in \IZ_{\geq 1}\times \IZ\times \IZ_{\geq 0}$ and fixed 
$\delta\in \IR_{>0}$ there is an algebraic moduli stack of finite type 
$\mfm_\delta^{ss}(\CX,r,e,v)$ 
of $\delta$-semistable objects of type $(r,e,v)$ of $\CC_\CX$. 
If $\delta<\delta'$ are two stability parameters so that 
$[\delta,\ \delta']\cap \Delta(r,e,v)=\emptyset$, the corresponding moduli stacks 
are canonically isomorphic. 
Moreover, for any $\delta\in \IR$ there are canonical open embeddings 
\be\label{eq:opembd}
\mfm_\delta^{ss}(\CX,r,e,v)\hookrightarrow {\mathfrak {Ob}}(\CX,r,e,v)
\hookrightarrow {\mathfrak {Ob}}(\CX).
\ee

\end{prop}

\subsection{ADHM invariants}\label{invariants} 
ADHM invariants will be defined applying the formalism of Joyce and 
Song \cite{genDTI} to $\delta$-semistable ADHM sheaves on $X$. Given corollary 
(\ref{extcoro}), the required results on Behrend constructible 
functions are a straightforward generalization of the analogous 
statements proven in 
  \cite[Sect. 7]{chamberI} for ADHM sheaves with $v=1$. 
Therefore the construction of generalized Donaldson-Thomas invariants 
via Behrend's constructible functions \cite{genDTI} applies to the present case. 

Let ${\sf L}(\CX)$ be the Lie algebra over $\IQ$ spanned by 
$\{\lambda(\gamma)\,|\, \gamma\in \IZ^3\}$ with Lie bracket 
\[
[\lambda(\gamma'), \lambda(\gamma'')] = (-1)^{\chi(\gamma',\gamma'')} 
\chi(\gamma',\gamma'') \lambda(\gamma'+\gamma'') 
\]
where 
\[
\chi(\gamma',\gamma'') = v''e'-v'e''-(v''r'-v'r'')(g-1)
\]
for any $\gamma'=(r',e',v')$, $\gamma''=(r'',e'',v'')$. 
Then there is a Lie algebra morphism 
\be\label{eq:Liemorphism}
\Psi: {\sf {SF}}^{\sf {ind}}_{\sf{al}}(\obj(\CX)) \to {\sf L}(\CX) 
\ee 
so that for any stack function of the form $[({\mathfrak X}, \rho)]$, whith
$\rho:{\mathfrak X} \hookrightarrow \obj(\CX,\gamma)\hookrightarrow \obj(\CX)$ 
an open embedding, 
and ${\mathfrak X}$ a $\IC^\times$-gerbe over an algebraic space 
${\sf X}$, 
\[
\Psi([({\mathfrak X}, \rho)]) = -\chi^B({\sf X}, \rho^*\nu) \lambda(\gamma) 
\]
where $\nu$ is Behrend's constructible function of the stack $\obj(\CX)$.

In order to define ADHM invariants note that for any $\delta\in \IR$, the canonical 
open embedding 
stack $\mfm_\delta^{ss}(\CX,\gamma)\hookrightarrow \obj(\CX)$ determines a stack 
function $\mfd_\delta(\gamma)\in {\underline {\sf{SF}}}(\obj(\CX))$. 
For $v=0$, the resulting stack functions are independent of stability 
parameters and will be denoted by $\mfh(\gamma)$. 

According to \cite[Thm. 8.7]{J-III} the associated log stack function 
\be\label{eq:logstfctB}
\mfe_\delta(\gamma) = \sum_{l\geq 1} {(-1)^{l-1}\over l} 
\mathop{\sum_{\gamma_1+\cdots+\gamma_l=\gamma}}_{\mu_\delta(\gamma_i)=
\mu_\delta(\gamma),\ 1\leq i\leq l} 
\mfd_\delta(\gamma_1)\ast\cdots\ast  \mfd_\delta(\gamma_l)
\ee
belongs to ${\sf {SF}}^{\sf {ind}}_{\sf{al}}(\obj(\CX))$, and is supported in 
$\obj(\CX,\gamma)$. 
Note that for fixed $\gamma$ and $\delta$ the sum in the right hand side 
is finite, therefore there are no convergence issues in the present case. 

Then, for $\gamma\in \IZ_{\geq 1}\times\IZ\times 
\IZ_{\geq 0}$,  the $\delta$-ADHM invariant $A_\delta(\gamma)$ is defined by 
\be\label{eq:ADHMinv} 
\Psi(\mfe_\delta(\gamma)) = -A_\delta(\gamma) \lambda(\gamma).
\ee
Note that 
$\mfe_\delta(\gamma)$ is independent of $\delta$ if $v=0$. 
Then the  corresponding 
invariants will be denoted by $H(\gamma)$.

By analogy with \cite{genDTI}, define the invariants ${\overline A}_\delta(r,e,v)$ 
by the multicover formula
\be\label{eq:multicoverA}
{\overline A}_\delta(r,e,v) = \mathop{\sum_{m\geq 1}}_{m|r,\ m|e,\ m|v}
{1\over m^2} {\overline A}_\delta(r/m, e/m,v/m).
\ee
Conjecturally, ${\overline A}_\delta(r/m, e/m,v/m)$ are integral. Obviously, 
for $v=0$ the alternative notation ${\overline H}(r,e)$ will be used. 

\section{Wallcrossing formulas}\label{wallsection} 

\subsection{Stack function identities}\label{sfidentities} 
Let $\gamma=(r,e,v)\in \IZ_{\geq 1}\times \IZ\times \IZ_{\geq 1}$ and let $\delta_c>0$ be 
a critical stability parameter of type $\gamma$. Let $\delta_- <\delta_c$, 
$\delta_+>\delta_c$ be stability parameters as in lemma (\ref{HNfiltrations}). 
Recall that ${\mathcal {HN}}_\pm(r,v)$ denote the set of positive, respectively negative 
admissible configurations of type $(r,v)$ introduced in definition (\ref{admconfig}). 
For any $h\in \IZ_{\geq 2}$ let ${\mathcal {HN}}_{\pm}(\gamma,\delta_c,h)$ denote the set of ordered sequences of 
triples $\left(\gamma_i =(r_i,e_i,v_i)\in \IZ_{\geq 1}\times \IZ\times \IZ_{\geq 0} \right)_{1\leq i\leq h}$ so that $\left(\rho_i=(r_i,v_i)\right)_{1\leq i\leq h}
\in {\mathcal{HN}}_\pm(r,v)$,
\[
e_1+\cdots+e_h =e \qquad \mathrm{and}
\qquad {e_i+v_i\delta_c\over r_i} = {e+v\delta_c\over r} \qquad 
\mathrm{for \ all} \ 1\leq i\leq h.
\]
More generally, given $h\in \IZ_{\geq 2}$, for any $0\leq k \leq h-1$ let  
${\mathcal {HN}}_{+}(\gamma,\delta_c,h,k)$ denote the set of ordered sequences 
$\left(\gamma_i =(r_i,e_i,v_i)\in \IZ_{\geq 1}\times \IZ\times \IZ_{\geq 0} \right)_{1\leq i\leq h}$ so that 
\begin{itemize}
\item $\gamma_1+\cdots +\gamma_h=\gamma$, $v_{h-k+1}=\cdots = v_h = 0$,
$v_i>0$ for $1\leq i \leq h-k$ , and 
$${e_1+ v_1\delta_c \over r_1} = \cdots ={e_{h-k}+v_{h-k}\delta_c
\over r_{h-k}} = {e_{h-k+1}\over r_{h-k+1}} = \cdots = {e_h \over r_h}={e+v\delta_c\over r}$$
\item The sequence $\left(\rho_j= (r_{j}, v_{j})\right)_{1\leq j\leq h-k}$ 
belongs to ${\mathcal {HN}}_{+}\left(r-\sum_{i=1}^k r_i, v\right)$. 
\end{itemize} 
Similarly, for any $0\leq k \leq h-1$ let  
${\mathcal {HN}}_{-}(\gamma,\delta_c,h,k)$ denote the set of ordered sequences 
$\left(\gamma_i =(r_i,e_i,v_i)\in \IZ_{\geq 1}\times \IZ\times \IZ_{\geq 0} \right)_{1\leq i\leq h}$ so that 
\begin{itemize}
\item $\gamma_1+\cdots +\gamma_h=\gamma$, $v_1=\cdots = v_k = 0$,
$v_i>0$ for $k+1\leq i \leq h$ , and 
$${e_1\over r_1} = \cdots ={e_{k}\over r_{k}} = {e_{k+1}+v_{k+1}\delta_{c}\over r_{k+1}} = \cdots = {e_h+v_h\delta_c \over r_h}={e+v\delta_c \over r}$$
\item The sequence $\left(\rho_j= (r_{k+j}, v_{k+j})\right)_{1\leq j\leq h-k}$ 
belongs to ${\mathcal {HN}}_{-}\left(r-\sum_{i=1}^k r_i, v\right)$. 
\end{itemize} 
\begin{rema}\label{configremaB} 
$(i)$ Obviously, in both cases $v_i>0$ for all $1\leq i\leq h$ if $k=0$. 
Moreover, \[{\mathcal{HN}}_\pm(\gamma,\delta_c,h) = 
{\mathcal{HN}}_\pm(\gamma,\delta_c,h,0) \cup  {\mathcal{HN}}_\pm(\gamma,\delta_c,h,1).
\]
If $k=h-1$ the condition that the sequence $(\rho_j)_{1\leq j \leq h-k}$ belong to 
${\mathcal {HN}}_{\pm}\left(r-\sum_{i=1}^k r_i, v\right)$ is empty.

$(ii)$ For fixed $\gamma$ and $\delta_c>0$ it straightforward to check that the 
following set is finite 
\[
\bigcup_{h\geq 2} \bigcup_{0\leq k \leq h-1} {\mathcal {HN}}_{\pm}(\gamma,\delta_c,h,k),
\]
i.e. the set ${\mathcal {HN}}_{\pm}(\gamma,\delta_c,h,k)$ is nonempty
only for a finite set of pairs $(h,k)$.
\end{rema}

For any triple $\gamma'=(r',e',v')\in \IZ_{\geq 1}\times \IZ\times \IZ_{\geq 1}$ 
let $\mfd_\pm(\gamma'),\mfd_{c}(\gamma')$ be the stack functions determined by the open embeddings 
$\mfm_{\delta_\pm}^{ss}(\CX,r',e',v')\hookrightarrow {\mathfrak {Ob}}(\CX)$, 
respectively $\mfm_{\delta_c}^{ss}(\CX,r',e',v')\hookrightarrow {\mathfrak {Ob}}(\CX)$.
The alternative notation $\mfh(\gamma')$ will be used if $v'=0$. 
\begin{lemm}\label{wallcrossinglemmaA} 
The following relations hold in the stack function algebra 
${\underline {\sf{SF}}}({\mathfrak {Ob}}(\CX))$ 
\be\label{eq:wallformulaA} 
\begin{aligned} 
\mfd_c(\gamma) = \mfd_{\pm}(\gamma) + \mathop{\sum_{h\geq 2}}_{}\
\mathop{\sum_{(\gamma_i) \in {\mathcal {HN}_\pm(\gamma, \delta_c,h)}}}_{}
\mfd_{\pm}(\gamma_1)\ast\cdots \ast\mfd_{\pm}(\gamma_h)
\end{aligned}
\ee
\be\label{eq:wallformulaB} 
\begin{aligned} 
& \mfd_-(\gamma) + \mathop{\sum_{h\geq 2}} \ 
\mathop{\sum_{(\gamma_i)\in {\mathcal{HN}}_-(\gamma,\delta_c,h,0)}}_{}
\mfd_{-}(\gamma_1)\ast\cdots \ast\mfd_{-}(\gamma_h)
= \\
& \mfd_c(\gamma) + \mathop{\sum_{h\geq 2}}(-1)^{h-1} 
\mathop{\sum_{(\gamma_i)\in {\mathcal {HN}}(\gamma,\delta_c,h,h-1)}}_{}
\mfh(\gamma_1) \ast \cdots \ast\mfh(\gamma_{h-1}) \ast \mfd_c(\gamma_h)\\
\end{aligned}
\ee
\end{lemm}

{\it Proof.} Equation \eqref{eq:wallformulaA} follows directly from lemma (\ref{HNfiltrations}). 
In order to prove formula \eqref{eq:wallformulaB} it will be first proven by induction 
that the following 
formula holds for any $l\in \IZ_{\geq 1}$. 
\be\label{eq:recformA} 
\bal
& \mfd_-(\gamma) + \mathop{\sum_{h\geq 2}} \ 
\mathop{\sum_{(\gamma_i)\in {\mathcal{HN}}_-(\gamma,\delta_c,h,0)}}_{}
\mfd_{-}(\gamma_1)\ast\cdots \ast\mfd_{-}(\gamma_h) = \\
& \mfd_{c}(\gamma) + \mathop{\sum_{k=2}^l}_{}(-1)^{k-1} \ 
\mathop{\sum_{(\gamma_i)\in {\mathcal {HN}}_-(\gamma,\delta_c,k,k-1)}}_{}
\mfh(\gamma_1) \ast \cdots \ast\mfh(\gamma_{k-1}) \ast \mfd_c(\gamma_k)\\
& + (-1)^l \mathop{\sum_{h\geq l+1}}_{}  \mathop{\sum_{(\gamma_i)\in 
{\mathcal{HN}}_-(\gamma, \delta_c, h, l)}}_{} \mfh(\gamma_1)\ast
\cdots \ast \mfh(\gamma_l) \ast \mfd_-(\gamma_{l+1}) \ast \cdots \ast
\mfd_-(\gamma_h)
\\
\eal
\ee
First note that remark (\ref{configremaB}.$ii$) implies that 
all sums in equation \eqref{eq:recformA} are finite for any $l\geq 1$.

Next, if $l=1$, equation \eqref{eq:recformA} is equivalent to \eqref{eq:wallformulaA}. 
Suppose it holds for some $l\geq 1$. 
Then note that equation \eqref{eq:wallformulaA} is valid for any triple $\gamma=(r,e,v)$ 
and any stability parameter $\delta_c$. If $\delta_c$ is not critical of type $\gamma$ 
as assumed above, it reduces to a trivial identity. 
In particular setting $\gamma= \gamma_{l+1}$ in 
equation \eqref{eq:wallformulaA} yields 
\[
\bal 
\mfd_{-}(\gamma_{l+1}) = 
& \mfd_{c}(\gamma_{l+1}) - \mathop{\sum_{m\geq 2}}_{} \ 
\mathop{\sum_{(\eta_i)\in 
{\mathcal {HN}}_-(\gamma_{l+1}, \delta_c, m, 1)}}_{} \mfh(\eta_1) \ast 
\mfd_-({\eta_2})\ast\cdots\ast\mfd_-({\eta_m})  \\   & - 
 \mathop{\sum_{m\geq 2}}_{} \ 
\mathop{\sum_{(\eta_i)\in 
{\mathcal {HN}}_-(\gamma_{l+1}, \delta_c, m, 0)}}_{} \mfd_-(\eta_1) \ast 
\mfd_-({\eta_2})\ast\cdots\ast\mfd_-({\eta_m}) 
\eal
\]
Using this expression, the third term in the right hand side of equation \eqref{eq:recformA}
can be rewritten as follows. 
\be\label{eq:recformB}
\bal 
& (-1)^l \mathop{\sum_{h\geq l+1}}_{}  \mathop{\sum_{(\gamma_i)\in 
{\mathcal{HN}}_-(\gamma, \delta_c, h, l)}}_{} \mfh(\gamma_1)\ast
\cdots \ast \mfh(\gamma_l) \ast \mfd_-(\gamma_{l+1}) \ast \cdots \ast
\mfd_-(\gamma_h) =\\
\eal
\ee
\[
\bal
&  (-1)^{l}
\mathop{\sum_{(\gamma_i)\in {\mathcal{HN}}_-(\gamma, \delta_c, l+1, l)}}_{}
\bigg[ 
\mfh(\gamma_1)\ast\cdots\ast\mfh(\gamma_{l}) \ast\mfd_c(\gamma_{l+1}) -  \\
& \mathop{\sum_{m\geq 2}}_{}\ \mathop{\sum_{(\eta_i)\in 
{\mathcal {HN}}_-(\gamma_{l+1}, \delta_c, m, 1)}}_{} 
\mfh(\gamma_1)\ast\cdots\ast\mfh(\gamma_{l})\ast \mfh(\eta_1)\ast 
\mfd_-({\eta_2})\ast\cdots\ast\mfd_-({\eta_m}) \\
& -\mathop{\sum_{m\geq 2}}_{}\ \mathop{\sum_{(\eta_i)\in 
{\mathcal {HN}}_-(\gamma_{l+1}, \delta_c, m, 0)}}_{} 
\mfh(\gamma_1)\ast\cdots\ast\mfh(\gamma_{l})\ast \mfd_-(\eta_1)\ast 
\mfd_-({\eta_2})\ast\cdots\ast\mfd_-({\eta_m}) \bigg]\\
\eal
\]
\[
+ (-1)^l \mathop{\sum_{h\geq l+2}}_{} \mathop{\sum_{(\gamma_i)\in 
{\mathcal{HN}}_-(\gamma, \delta_c, h, l)}}_{} \mfh(\gamma_1)\ast
\cdots \ast \mfh(\gamma_l) \ast \mfd_-(\gamma_{l+1}) \ast \cdots \ast
\mfd_-(\gamma_h)
\]
By construction 
\[
\mathop{{\bigcup}_{(\gamma_i)\in {\mathcal{HN}}_-(\gamma, \delta_c, l+1, l)}}_{}
{\mathcal {HN}}_-(\gamma_{l+1},\delta_c, m, j) = 
{\mathcal{HN}}_-(\gamma, \delta_c, l+m, l+j)
\]
for any $m\in \IZ_{\geq 2}$, $j\in \{0,1\}$. Therefore the last two terms in the right hand side of equation \eqref{eq:recformB} cancel, and formula \eqref{eq:recformB}
reduces to 
\be\label{eq:recformC}
\bal 
& (-1)^l \mathop{\sum_{h\geq l+1}}_{}  \mathop{\sum_{(\gamma_i)\in 
{\mathcal{HN}}_-(\gamma, \delta_c, h, l)}}_{} \mfh(\gamma_1)\ast
\cdots \ast \mfh(\gamma_l) \ast \mfd_-(\gamma_{l+1}) \ast \cdots \ast
\mfd_-(\gamma_h) =\\
& (-1)^{l}
\mathop{\sum_{(\gamma_i)\in {\mathcal{HN}}_-(\gamma, \delta_c, l+1, l)}}_{}
\mfh(\gamma_1)\ast\cdots\ast\mfh(\gamma_{l}) \ast\mfd_c(\gamma_{l+1}) -  \\
& +(-1)^{l+1} \mathop{\sum_{h\geq l+2}}_{} \mathop{\sum_{(\gamma_i)\in 
{\mathcal{HN}}_-(\gamma, \delta_c, h, l+1)}}_{} \mfh(\gamma_1)\ast
\cdots \ast \mfh(\gamma_{l+1}) \ast \mfd_-(\gamma_{l+2}) \ast \cdots \ast
\mfd_-(\gamma_h)
\eal
\ee
Substituting \eqref{eq:recformC} in \eqref{eq:recformA} it follows that formula 
\eqref{eq:recformA} also holds if $l$ is replaced by $(l+1)$. This concludes the inductive 
proof of formula \eqref{eq:recformA}. 

In order to conclude the proof of equation \eqref{eq:wallformulaB}, it suffices to observe 
that for sufficiently large $l$, equation \eqref{eq:recformA} stabilizes to equation 
\eqref{eq:wallformulaB} using remark (\ref{configremaB}.$ii$). 

\hfill $\Box$

Now note that equations \eqref{eq:wallformulaA}, \eqref{eq:wallformulaB} yield a 
recursive algorithm expressing $\mfd_-(\gamma)$ in terms of 
$\mfd_+(\gamma_i)$, $1\leq i\leq h$, $h\geq 1$. This follows observing 
that in the left hand side of \eqref{eq:wallformulaB} $0<v_i<v$ for all 
stack functions $\mfd_-(\gamma_i)$ occuring in the sum 
\[
\mathop{\sum_{h\geq 2}} \ 
\mathop{\sum_{(\gamma_i)\in {\mathcal{HN}}_-(\gamma,\delta_c,h,0)}}_{}
\mfd_{-}(\gamma_1)\ast\cdots \ast\mfd_{-}(\gamma_h).
\]
Therefore, once a formula for the difference $\mfd_-(\gamma)-\mfd_+(\gamma)$, 
has been derived for triples of the form 
$\gamma=(r,e,v)$, one can recursively derive an analogous formula for triples of the form 
$\gamma=(r,e,v+1)$. 
For $v=1$, equations \eqref{eq:wallformulaA}, \eqref{eq:wallformulaB} 
easily imply 
\be\label{eq:wallformulaD}
\bal 
 \mfd_-(\gamma) =  \mfd_+(\gamma) & +\mathop{\sum_{l\geq 2}}_{}
 (-1)^{l} \mathop{\sum_{(\gamma_i)\in {\mathcal{HN}}_-(\gamma, 
\delta_c,l,l-1)}}_{} \mfh(\gamma_1)\ast\cdots\ast[\mfd_+(\gamma_l),\mfh(\gamma_{l-1})]\\
 \eal
 \ee 
Employing the above recursive algorithm one can determine in principle 
analogous formulas for $v\geq 2$. 
Since the resulting expressions quickly become cumbersome, 
explicit formulas will be given below only for $v=2$. 

\begin{coro}\label{ranktwocor} 
Suppose $\gamma=(r,e,2)$ with $(r,e)\in \IZ_{\geq 1}\times \IZ$. 
The following relations hold in the stack function algebra 
${\underline {\sf{SF}}}({\mathfrak {Ob}}(\CX))$ 
\be\label{eq:wallformulaC}
\bal
& \mfd_-(\gamma) = \mfd_+(\gamma) + \mathop{\sum_{l\geq 2}}_{}
 (-1)^{l} \mathop{\sum_{(\gamma_i)\in {\mathcal{HN}}_-(\gamma, 
\delta_c,l,l-1)}}_{} \mfh(\gamma_1)\ast\cdots\ast[ \mfd_+(\gamma_{l}),
\mfh(\gamma_{l-1}) ]\\
& 
+ \mathop{ \sum_{(\gamma_1,\gamma_2)\in {\mathcal{HN}}_{+} 
(\gamma, \delta_c, 2,0)}}_{} \mfd_+(\gamma_1)\ast\mfd_+(\gamma_2)
-\mathop{\sum_{(\gamma_1,\gamma_2)\in {\mathcal{HN}}_-(\gamma, 
\delta_c,2,0)}}_{} \mfd_-(\gamma_1) \ast \mfd_-(\gamma_2) \\
& + \mathop{\sum_{l\geq 2}}_{}
 (-1)^{l} \mathop{\sum_{(\gamma_i)\in {\mathcal{HN}}_-(\gamma, 
\delta_c,l+1,l-1)}}_{} 
\mfh(\gamma_1)\ast\cdots\ast[ \mfd_+(\gamma_{l+1})\ast \mfd_+(\gamma_{l}),
\mfh(\gamma_{l-1}) ]
\eal
\ee
where $\mfd_-(\gamma_1), \mfd_-(\gamma_2)$ are given by 
equation \eqref{eq:wallformulaD}. 
\end{coro}

\subsection{Wallcrossing for $v=2$ invariants}\label{wallformula}
Let $\gamma=(r,e,2)$, $(r,e)\in \IZ_{\geq 1}\times \IZ$, 
$\delta_c>0$ a critical stability parameter of type $\gamma$, and 
$\delta_\pm$ two noncritical stability parameters 
as in lemma (\ref{HNfiltrations}). The main goal of this section is 
to convert the stack function relation \eqref{eq:wallformulaC} to a
wallcrossing formula for generalized Donaldson-Thomas invariants 
of ADHM sheaves. 

As mentioned in the introduction the alternative notation $\alpha =(r,e)$ will be used 
for 
pairs $(r,e)\in \IZ_{\geq 1}\times \IZ$. Using this notation, the sets 
${\mathcal{HN}}_-(\alpha,v, \delta_c, h,k)$, $v\in\{1,2\}$, $k\in\{0, h-2,h-1\}$, 
can be identified with sets of ordered sequences 
$(\alpha_i)_{1\leq i\leq h}$ satisfying the conditions listed above theorem 
(\ref{wallcrossingthmA}). For convenience, recall that ${\mathcal {HN}}_-(\alpha,v,\delta_c,l,l-1)$, $l\in \IZ_{\geq 1}$, $v\in \{1,2\}$, denotes 
the set of ordered sequences
$((\alpha_i))_{1\leq i\leq l}$, $\alpha_i\in \IZ_{\geq 1}\times \IZ$,
$1\leq i\leq l$ so that 
\be\label{eq:alphadecompX}
\alpha_1+\cdots +\alpha_l=\alpha
\ee
and 
\be\label{eq:slopesXA}
{e_1\over r_1} = \cdots = {e_{l-1}\over r_{l-1}}={e_l+v\delta_c\over r_l}
={e+v\delta_c\over r}
\ee
Similarly, ${\mathcal {HN}}_-(\alpha,v,\delta_c,l,l-2)$,
 $l\in \IZ_{\geq 2}$, 
 denotes the set of ordered sequences
$((\alpha_i))_{1\leq i\leq l}$, $\alpha_i\in \IZ_{\geq 1}\times \IZ$,
$1\leq i\leq l$ satisfying condition \eqref{eq:alphadecomp},
\be\label{eq:slopesB}
{e_1\over r_1} = \cdots = {e_{l-2}\over r_{l-2}}={e_{l-1}+\delta_c\over r_{l-1}}
={e_{l}+\delta_c\over r_l} = {e+2\delta_c\over r} 
\ee
and $1/{r_{l-1}} < 1/{r_{l}}$. 

Note that the sets ${\mathcal{HN}}_-(\alpha,2, \delta_c, h,0)$ are nonempty if and only if 
$h=2$, in which case they consist of ordered pairs $(\alpha_1,\alpha_2)$ so that 
$\alpha_1+\alpha_2=\alpha$, $1/r_1 < 1/r_2$, and 
\[
{e_1+\delta_c\over r_1} ={e_2+\delta_c\over r_2} = {e+2\delta_c\over r}
\]
Moreover the set ${\mathcal{HN}}_-(\alpha,2, \delta_c, 1,0)$ consists of 
only of the element $(\alpha)$. 
 
It is straightforward to check that for fixed $\alpha=(r,e)$ and $\delta_c$, 
the union 
\be\label{eq:unionset}
\bal 
& \qquad\qquad \bigcup_{l\geq 1}
 \left[{\mathcal {HN}}_-(\alpha,2,\delta_c,l,l-1) \cup  
 {\mathcal {HN}}_-(\alpha,2,\delta_c,l+1,l-1)\right]\\
 &\bigcup_{(\alpha_1,\alpha_2)\in 
 {\mathcal {HN}}_-(\alpha,2,\delta_c,2,0) }
 \bigcup_{l_1\geq 1}\bigcup_{l_2\geq 1}
 \left[{\mathcal {HN}}_-(\alpha_1,1,\delta_c,l_1,l_1-1) \times 
 {\mathcal {HN}}_-(\alpha_2,1,\delta_c,l_2,l_2-1) \right]
  \\
  \eal
 \ee
 is a finite set.

Now let $0<\delta_-<\delta_c<\delta_+$ be stability parameters so that there are no 
critical stability parameters of type $(\alpha,2)$ in the intervals $[\delta_-,\ \delta_c)$, 
$(\delta_c,\ \delta_+]$. Since the set \eqref{eq:unionset} is finite $\delta_-, \delta_+$ 
can be chosen so that the same holds for all numerical types $(\alpha_i, v_i)$ 
in all ordered sequences in \eqref{eq:unionset}. 
Then the following lemma holds. 
\begin{lemm}\label{wallcrossinglemmaB}
The following relations hold in the stack function algebra ${\underline{\sf {SF}}}(\obj(\CX))$
\be\label{eq:wallformulaE}
\bal 
& \mfd_-(\alpha,1)= \mathop{\sum_{l\geq 1}}_{} {(-1)^{l-1}\over (l-1)!}
\mathop{\sum_{(\alpha_i)\in  {\mathcal{HN}}_-(\alpha,1, \delta_c, l,l-1)}}_{}
[\mfg(\alpha_1),[\cdots[\mfg(\alpha_{l-1}),\mfd_+(\alpha_l,1)]\cdots]\\
\eal 
\ee
\be\label{eq:wallformulaF}
\bal 
&\mfd_-(\alpha,2) = \mathop{\sum_{l\geq 1}}_{} {(-1)^{l-1}\over (l-1)!}
\mathop{\sum_{(\alpha_i)\in  {\mathcal{HN}}_-(\alpha,2, \delta_c, l,l-1)}}_{}
[\mfg(\alpha_1),[\cdots[\mfg(\alpha_{l-1}),\mfd_+(\alpha_l,2)]\cdots]\\
& +\mathop{\sum_{l\geq 1}}_{} {(-1)^{l-1}\over (l-1)!}
\mathop{\sum_{(\alpha_i)\in  {\mathcal{HN}}_-(\alpha,2, \delta_c, l+1,l-1)}}_{}
[\mfg(\alpha_1),[\cdots[\mfg(\alpha_{l-1}),\mfd_+(\alpha_{l+1},1)\ast
\mfd_+(\alpha_{l},1)]\cdots]\\
& -\mathop{\sum_{(\alpha_1,\alpha_2)\in 
{\mathcal{HN}}_-(\alpha,2,\delta_c,2,0)}}_{} \mathop{\sum_{l_1\geq 1}}_{}
\mathop{\sum_{l_2\geq 1}}_{} {(-1)^{l_1-1}\over (l_1-1)!}{(-1)^{l_2-1}\over (l_2-1)!}
\mathop{\sum_{(\alpha_{1,i})\in  {\mathcal{HN}}_-(\alpha_1,1, \delta_c, l_1,l_1-1)}}_{}\\
& 
\mathop{\sum_{(\alpha_{2,i})\in  {\mathcal{HN}}_-(\alpha_2,1, \delta_c, l_2,l_2-1)}}_{}
 \big(
 [\mfg(\alpha_{1,1}),[\cdots[\mfg(\alpha_{1,l_1-1}),\mfd_+(\alpha_{1,l_1},1)]\cdots]\\
 & \qquad \qquad\qquad \qquad \qquad \ 
 \ast[\mfg(\alpha_{2,1}),[\cdots[\mfg(\alpha_{2,l_2-1}),\mfd_+(\alpha_{2,l_2},1)]\cdots]
 \big).\\
\eal
\ee
\end{lemm}

{\it Proof.} Formulas \eqref{eq:wallformulaE}, 
\eqref{eq:wallformulaF} follow from equations  \eqref{eq:wallformulaC}, 
\eqref{eq:wallformulaD} by repeating the computations 
in the proof of  \cite[Lemm. 2.6]{chamberI} in the present context. 

\hfill $\Box$

{\it Proof of Theorem (\ref{wallcrossingthmA}.)} 
The proof consists of two steps. First the stack function identities \eqref{eq:wallformulaE}, 
\eqref{eq:wallformulaF} must be converted into similar identities for the log 
stack functions \eqref{eq:logstfctB}. As explained in \cite[Sect. 6.5]{J-IV}, 
\cite[Sect. 3.5]{genDTI}, 
applying the morphism 
\eqref{eq:Liemorphism} to the log stack function identities \eqref{eq:wallformulaE}, 
\eqref{eq:wallformulaF} yields certain relations 
in the universal enveloping algebra $U({\sf L}(\CX))$ of the Lie algebra ${\sf L}(\CX)$. 
These relations imply in turn a wallcrossing formula for generalized Donaldson-Thomas 
invariants by identifying the coefficients of generators of the generators of 
${\sf L}(\CX)\subset U({\sf L}(\CX))$. 

Given the above choice of $\delta_\pm$, for 
 $v=1$, equation \eqref{eq:logstfctB} reduces to $\mfe_\pm(\gamma)=
\mfd_\pm(\gamma)$, while for $v=2$ 
\be\label{eq:logstfctA} 
\bal 
\mfe_\pm(\gamma) = \mfd_\pm(\gamma) -{1\over 2} \mfd_\pm(\gamma/2)\ast \mfd_\pm(\gamma/2).
\eal
\ee
The second term in the right hand side of \eqref{eq:logstfctA}
is by convention trivial unless $(r,e)$ are even. 

Equations \eqref{eq:logstfctA}, \eqref{eq:wallformulaF}, \eqref{eq:wallformulaG}
yield the following identity in the universal enveloping algebra of the Lie 
algebra ${\sf L}(\CX)$ 
\be\label{eq:univid}
\bal 
& \mathop{\sum_{\alpha}}(A_{-}(\alpha,2)-A_{+}(\alpha,2))\lambda(\alpha,2) = \\
&\mathop{\sum_{\alpha}}\mathop{\sum_{l\geq 2}}_{} {1\over (l-1)!}
\mathop{\sum_{(\alpha_i)\in  {\mathcal{HN}}_-(\alpha,2, \delta_c, l,l-1)}}_{}
\left(A_+(\alpha_l,2)\prod_{i=1}^{l-1}f_2(\alpha_i)H(\alpha_i)\right) \lambda(\alpha,2)\\
\eal 
\ee
\[
\bal
&-\mathop{\sum_{\alpha}}\mathop{\sum_{l\geq 1}}_{} {1\over (l-1)!}
\mathop{\sum_{(\alpha_i)\in  {\mathcal{HN}}_-(\alpha,2, \delta_c, l+1,l-1)}}_{}
\left(A_+(\alpha_l,1)
A_+(\alpha_{l+1},1)
\prod_{i=1}^{l-1}H(\alpha_i)\right)\\
& \qquad \qquad \qquad \qquad\qquad \qquad \quad\
[\lambda(\alpha_1),[\cdots[\lambda(\alpha_{l-1}), 
\lambda(\alpha_{l+1},1)\star\lambda(\alpha_{l},1)]\cdots ]\\
& +\mathop{\sum_{\alpha}}\mathop{\sum_{(\alpha_1,\alpha_2)\in 
{\mathcal{HN}}_-(\alpha,2,\delta_c,2,0)}}_{} \mathop{\sum_{l_1\geq 1}}_{}
\mathop{\sum_{l_2\geq 1}}_{} {1\over (l_1-1)!}{1\over (l_2-1)!}
\mathop{\sum_{(\alpha_{1,i})\in  {\mathcal{HN}}_-(\alpha_1,1, \delta_c, l_1,l_1-1)}}_{}\\
& \mathop{\sum_{(\alpha_{2,i})\in  {\mathcal{HN}}_-(\alpha_2,1, \delta_c, l_2,l_2-1)}}_{}
A_+(\alpha_{1,l_1})A_+(\alpha_{2,l_2}) 
\prod_{i=1}^{l_1-1} f_1(\alpha_{1,i})H(\alpha_{1,i}) 
\prod_{i=1}^{l_2-1} f_1(\alpha_{2,i})H(\alpha_{2,i}) \\
& \qquad\qquad\qquad\qquad\qquad\quad
\lambda(\alpha_1,1)\star\lambda(\alpha_2,1)\\
&+{1\over 2}\mathop{\sum_{\alpha}} (A_{-}(\alpha/2,1)^2-A_{+}(\alpha/2,1)^2) 
\lambda(\alpha/2,1)\star \lambda(\alpha/2,1) \\
&-\mathop{\sum_{\alpha}}\mathop{\sum_{l\geq 1}}_{} {1\over (l-1)!}
\mathop{\sum_{(\alpha_i)\in  {\mathcal{HN}}_-(\alpha,2, \delta_c, l,l-1)}}_{}
\left(A_+(\alpha_l/2,1)^2
\prod_{i=1}^{l-1}H(\alpha_i)\right) \\
& \qquad\qquad\qquad\qquad
[\lambda(\alpha_1),[\cdots,[\lambda(\alpha_{l-1}, \lambda(\alpha_l/2,1)\star 
\lambda(\alpha_l/2,1)]\cdots]\\
\eal
\]
where $\star$ denotes the associative product in the universal enveloping algebra. 
By conventions the invariants of the form $A_+(\alpha/2,1)$ are trivial unless 
$\alpha=2\alpha'$ for some $\alpha'=(r',e')=\IZ_{\geq 1} \times \IZ$. 

Next, the 
identity \cite[Eqn. 127]{J-IV} or \cite[Eqn. 45]{genDTI} yields the following relations in the universal enveloping 
algebra 
\[
\bal 
\lambda(\alpha_{l+1},1)\star\lambda(\alpha_{l},1) & = {1\over 2}
g(\alpha_{l+1},\alpha_{l}) 
\lambda(\alpha_l+\alpha_{l+1},2) + \cdots  \\
\lambda(\alpha_1,1)\star\lambda(\alpha_{2},1) & = {1\over 2}g(\alpha_1,\alpha_{2}) 
\lambda(\alpha_1+\alpha_2,2) + \cdots \\
\lambda(\alpha/2,1)\star\lambda(\alpha/2,1) & = \cdots\\
\lambda(\alpha_l/2,1)\star\lambda(\alpha_l/2,1) & = \cdots \\
\eal
\]
where $\cdots$ stands for linear combinations of generators of $U({\sf L}(\CX))$
not in ${\sf L}(\CX)$. 
Since the left hand side of equation \eqref{eq:univid} must belong to the Lie algebra 
${\sf L}(\CX)$ according to \cite[Thm. 8.7]{J-III}, it follows that all higher order terms 
must cancel. Then equation \eqref{eq:wallformulaG} follows by 
straightforward computations. 

\hfill $\Box$

\section{Comparison with Kontsevich-Soibelman Formula}\label{KSsect}

The goal of this section is to prove that 
formula \eqref{eq:wallformulaG} is in agreement with  the wallcrossing 
formula of 
Kontsevich and Soibelman \cite{wallcrossing}, which will be referred to as 
the KS formula in the following. 

As in section (\ref{wallformula}), numerical types of ADHM sheaves
will be denoted by $\gamma=(\alpha,v)$, $\alpha=
(r,e)\in \IZ_{\geq 1}\times \IZ$, $v\in \IZ_{\geq 0}$.
In order to streamline the computations, let 
${\sf L}(\CX)_{\leq 2}$ denote the truncation of the Lie algebra ${\sf L}(\CX)$ 
defined by 
\be\label{eq:chargeliebracket}
\begin{aligned}
& [\lambda(\alpha_1,v_1), \lambda(\alpha_2,v_2)]_{\leq 2} = 
\left\{\begin{array}{ll} 
[\lambda(\alpha_1,v_1), \lambda(\alpha_2,v_2)] & \quad \mathrm{if}\ v_1+v_2\leq 2\\
0 & \quad \mathrm{otherwise}.\\
\end{array}\right. \\
\end{aligned}
\ee
Furthermore, it will be more convenient to use the alternative notation 
${\sf e}_\alpha=\lambda(\alpha,0)$, ${\sf f}_\alpha=\lambda(\alpha,1)$,
 and ${\sf g}_{\alpha}=\lambda(\alpha,2)$. 
 
Given a 
critical stability parameter $\delta_c$ of type $(r,e,2)$, $(r,e)\in \IZ_{\geq 1}\times \IZ$, 
there exist two pairs
 $\alpha =(r_{\alpha}, e_{\alpha})$ and $\beta =(r_{\beta}, e_{\beta})$
with 
\[
\frac{e_{\alpha} + \delta_c}{r_{\alpha}}=\frac{e_{\beta}}{r_{\beta}}= \mu_{\delta_c}(\gamma)
\]
so that any $\eta\in \IZ_{\geq 1}\times \IZ$ 
with $\mu_{\delta_c}(\eta) = \mu_{\delta_c}(\gamma)$ can 
be uniquely written as $\eta = (q\beta,0), (\alpha+ q\beta,1),$ or 
$(2\alpha+ q\beta,2)$, with $q\in \IZ_{\geq 0}$.

For any $q\in \IZ_{\geq 0}$ the following formal expressions will 
be needed in the KS formula,
\begin{equation}
U_{\alpha+q\beta} = \mathrm{exp}({\sf f}_{\alpha+q\beta} + \frac{1}{4} {\sf g}_{2\alpha+2q\beta}) \ \ , \ \ 
U_{2\alpha+q\beta} = \mathrm{exp}({\sf g}_{2\alpha+q\beta}) \ \ , \ \
U_{q\beta} = \text{exp} (\sum_{m\geq1} \frac{{\sf e}_{mq\beta}}{m^2})\ .
\end{equation}
Moreover, let
\[
\IH = \sum_{q\geq 0} { H}(q\beta){{\sf e}_{q\beta}},
\]
where the invariants $H(\alpha)$ are defined in \eqref{eq:ADHMinv}.
Then the wallcrossing formula of Kontsevich and Soibelman
reads 
\begin{equation}\label{eq:KSformA}
\begin{aligned}
& \mathrm{exp(\IH)}\prod_{{q\geq 0},\ q\downarrow} U_{2\alpha + q\beta}^{\overline{A}_{+}(2\alpha+q\beta,2)} 
\prod_{{q\geq 0},\ q\downarrow} U_{\alpha+q\beta}^{A_{+}(\alpha+q\beta,1)}  \\
&= \prod_{q\geq 0,\ q\uparrow} U_{\alpha+ q\beta}^{A_{-}(\alpha+q\beta,1)}
\prod_{{q\geq 0},\ q\uparrow} U_{2\alpha + q\beta}^{\overline{A}_{-}(2\alpha+q\beta,2)} 
\mathrm{exp(\IH)}\end{aligned}
\end{equation}
where an up, respectively down arrow means that the factors in the corresponding 
product are taken in increasing, respectively decreasing order of $q$. Note that
$\overline{A}_{\pm}(2\alpha+q\beta,2)$ are the 
invariants defined in section (\ref{invariants})
by the multicover formula \eqref{eq:multicoverA}. In this case equation 
\eqref{eq:multicoverA} reduces to 
\[
 A_{\pm}(2\alpha+q\beta,2) = \overline{A}_{\pm}(2\alpha+q\beta,2) + \frac{1}{4} A_{\pm}(\alpha+q\beta/2,1).
\]
Expanding the right hand side, equation 
\eqref{eq:KSformA} yields 
\begin{equation}\label{eq:KSformB}
\begin{aligned}
& \mathrm{exp}(\sum_{q\geq0} A_{-}(2\alpha+q\beta,2) {\sf g}_{2\alpha+q\beta} + \\ 
& \sum_{q_2>q_1\geq0}
\frac{1}{2} g(q_1\beta,q_2\beta)A_{-}(\alpha+q_1\beta,1)A_{-}(\alpha+q_2\beta,1){\sf g}_{2\alpha+(q_1+q_2)\beta} )= \\
&  \mathrm{exp(\IH)}\, \mathrm{exp}(\sum_{q\geq0} A_{+}(2\alpha+q\beta,2) {\sf g}_{2\alpha+q\beta} \\
& + \sum_{q_1>q_2\geq0} \frac{1}{2} g(q_1\beta,q_2\beta)A_{+}(\alpha+q_1\beta,1)A_{+}(\alpha+q_2\beta,1) {\sf g}_{2\alpha+(q_1+q_2)\beta})\, \mathrm{exp(-\IH)},
\end{aligned}
\end{equation}
modulo terms involving ${\sf f}_{\gamma}$. These terms are omitted since they enter 
$v=1$ wallcrossing formula derived in 
\cite{chamberII}.
The BCH formula 
\begin{equation}
\begin{aligned}
\text{exp}(A) \text{exp}(B) \text{exp} (-A) & = \text{exp}( \sum_{n=0} \frac{1}{n!} (Ad(A))^n B )\\ & 
= \text{exp} ( B + [A,B] +\frac{1}{2} [A,[A,B]]+ \cdots),\\
\end{aligned}
\end{equation}
 yields 
\be\label{eq:BCHg}
\begin{aligned}
& \mathrm{exp}(\IH)\, \mathrm{exp}({\sf g}_{2\alpha+q\beta})\, \mathrm{exp}(-\IH) = 
\mathrm{exp}( {\sf g}_{2\alpha+q\beta} + 
\sum_{q_1>0} f_2(q_1\beta)H(q_1\beta) {\sf g}_{2\alpha+(q+q_1)\beta} \\
& + \frac{1}{2 !}\sum_{q_1>0,q_2>0} f_2(q_1\beta)H(q_1\beta) f_2(q_2\beta)H(q_2 \beta)  {\sf g}_{2\alpha+(q+q_1+q_2)\beta} +\cdots ) \\
& = \mathrm{exp} \Big( \sum_{l\geq0, q_i>0} \frac{1}{l !} (\prod_{i=1}^{l} f_2(q_i \beta) H(q_i \beta) ) 
{\sf g}_{2\alpha+(q+q_1+\cdots+q_l)\beta} \Big)
\end{aligned} 
\ee 
Substituting \eqref{eq:BCHg} in \eqref{eq:KSformB} 
results in 
\be\label{KSformC}
\begin{aligned}
&\mathrm{exp}\big( \sum_{q\geq0} A_{-}(2\alpha+q\beta,2){\sf g}_{2\alpha+q\beta} 
+\sum_{q_2>q_1\geq0} \frac{1}{2} g(q_1\beta, q_2\beta) A_{-}(\alpha+q_1\beta,1)
A_{-}(\alpha+q_2\beta,1) {\sf g}_{2\alpha+(q_1+q_2)\beta} \Big) \\
&= \mathrm{exp}\big( \sum_{ \substack{q\geq0, l\geq 0 \\ q_i>0}} A_{+}(2\alpha+q\beta,2) 
\frac{1}{l !} (\prod_{i=1}^{l} f_2(q_i \beta) H(q_i \beta) ) {\sf g}_{2\alpha+(q+q_1+\cdots+q_l)\beta} \\
& + \sum_{\substack{q_1' > q_2'\geq 0\\ l\geq0, q_i>0}} \frac{1}{2} g(q_1'\beta, q_2'\beta) 
A_{+}(\alpha+q_1'\beta,1) A_{+}(\alpha+q_2'\beta,1)  \frac{1}{l !}(\prod_{i=1}^{l} f_2(q_i \beta) H(q_i \beta) )
{\sf g}_{2\alpha+(q_1'+q_2'+q_1+\cdots+q_l)\beta} \Big)
\end{aligned}
\ee
In order to further simplify the notation, let 
\[
A_{\pm} (v\alpha+q\beta,v) \equiv A_{\pm}(q,v), \qquad {\sf g}_{2\alpha+q\beta}
\equiv {\sf g}_q .
\]
Comparing the coefficients of ${\sf g}_Q$ in \eqref{eq:KSformB},
yields 
\be\label{eq:KSformC}
\begin{aligned}
& A_{-}(Q,2) = \sum_{\substack{q'\geq0 ,\ l \geq0,\ q_i>0 \\ q'+q_1+\cdots+q_l=Q}} A_{+}(q',2) 
\frac{1}{l !} (\prod_{i=1}^{l} f_2(q_i \beta) H(q_i \beta) )  \\
& +\frac{1}{2}\sum_{\substack{q_1' > q_2'\geq 0\\ l\geq0,\ q_i>0 \\ q_1'+q_2'+ q_1+\cdots + q_l=Q}} 
 g(q_1'\beta, q_2'\beta) A_{+}(q_1',1) A_{+}(q_2',1)\frac{1}{l !}(\prod_{i=1}^{l} f_2(q_i \beta) H(q_i \beta) ) \\
& - \frac{1}{2} \sum_{ q_2' > q_1'\geq 0,\ q_1'+q_2'=Q}\ g(q_1'\beta, q_2'\beta) A_{-}(q_1',1)A_{-}(q_2',1) \ .
\end{aligned}
\ee
Using the $v=1$ wallcrossing formula \cite[Thm. 1.1]{chamberII} the last term in 
\eqref{eq:KSformC} becomes 
\be\label{eq:lastterm}
\begin{aligned}
& - \frac{1}{2} \sum_{ q_2 > q_1\geq 0,\ q_1+q_2=Q}
\ g(q_1\beta, q_2\beta) A_{-}(q_1,1)A_{-}(q_2,1) \\
& = - \frac{1}{2} \sum_{ \substack{ q_2> q_1\geq 0\\ q_1+q_2=Q\\ l \geq 0,\ \tilde{l} \geq 0\\q_1' \geq 0,\ q_2' \geq 0\\ n_i>0,\ \tilde{n}_i >0 \\ q_1'+n_1+\cdots+n_l = q_1 \\
q_2'+\tilde{n}_1+\cdots+\tilde{n}_{\tilde{l}} = q_2 }} g(q_1\beta, q_2\beta) A_{+}(q_1',1)A_{+}(q_2',1)
\frac{1}{l !} (\prod_{i=1}^{l} f_1(n_i \beta) H(n_i \beta) ) \frac{1}{\tilde l !} 
(\prod_{i=1}^{\tilde l} f_1(\tilde{n}_i \beta) H(\tilde{n}_i \beta) ) \ .
\end{aligned}
\ee 
Therefore the final wallcrossing formula for $v=2$ invariants is 
\be\label{eq:KSformD}
\begin{aligned}
& A_{-}(Q,2) = \sum_{\substack{q'\geq0,\ l \geq0,\ q_i>0 \\ q'+q_1+\cdots+q_l=Q}} A_{+}(q',2) 
\frac{1}{l !} (\prod_{i=1}^{l} f_2(q_i \beta) H(q_i \beta) )  \\
& +\frac{1}{2}\sum_{\substack{q_1' > q_2'\geq 0\\ l\geq0,\ q_i>0 \\ q_1'+q_2'+ q_1+\cdots + q_l=Q}} 
\frac{1}{2} g(q_1'\beta, q_2'\beta) A_{+}(q_1',1) A_{+}(q_2',1)\frac{1}{l !} (\prod_{i=1}^{l} f_2(q_i \beta) H(q_i \beta) ) \\
&  -\frac{1}{2} \sum_{ \substack{ q_2 > q_1\geq 0\\ q_1+q_2=Q\\ l \geq 0,\ \tilde{l} \geq 0\\q_1' \geq 0,\ q_2' \geq 0\\ n_i>0, \tilde{n}_i >0 \\ q_1'+n_1+\cdots+n_l = q_1 \\
q_2'+\tilde{n}_1+\cdots+\tilde{n}_{\tilde{l}} = q_2 }} g(q_1\beta, q_2\beta) A_{+}(q_1',1)A_{+}(q_2',1)
\frac{1}{l !} (\prod_{i=1}^{l} f_1(n_i \beta) H(n_i \beta) ) \frac{1}{\tilde l !} 
(\prod_{i=1}^{\tilde l} f_1(\tilde{n}_i \beta) H(\tilde{n}_i \beta) ) \ .
\end{aligned}
\ee
This formula agrees with \eqref{eq:wallformulaG} since the bilinear function 
$g({\quad},{\quad})$ is antisymmetric.

\section{Asymptotic invariants in the $g=0$ theory}\label{genuszero} 
In this subsection $X$ will be a smooth genus 0 curve over a $\IC$-field $K$, 
and $M_1\simeq \CO_X(d_1)$, $M_2\simeq \CO_X(d_2)$, with 
$(d_1,d_2)=(1,1)$ or $(d_1,d_2)=(0,2)$. In this case any coherent 
locally free 
sheaf $E$ on $X$ is isomorphic to a direct sum of line bundles. 
Let $E_{\geq 0}$ denote the direct sum of all summands of non-negative 
degree, and $E_{<0}$ denote the direct sum of all summands of 
negative degree. 

\begin{lemm}\label{zerophi} 
Let $\CE=(E,V,\Phi_1,\Phi_2,\phi,\psi)$ be a nontrivial $\delta$-semistable 
ADHM sheaf of type 
$(r,e,v)\in \IZ_{\geq 1}\times \IZ\times \IZ_{\geq 1}$, for some $\delta> 0$.
Then $E_{<0}=0$ 
 and $\phi$ is identically zero. 
\end{lemm} 

{\it Proof.}
Since $\delta>0$, lemma (\ref{basicprop}.$ii$) implies that $\psi$ 
is not identically zero. 
Then obviously $E_{\geq 0}$ must be nontrivial and $\mathrm{Im}(\psi) 
\subseteq E_{\geq 0}$. Since $M\simeq K_X^{-1}\simeq \CO_X(2)$, 
$E_{\geq 0}\otimes_X M\subseteq \mathrm{Ker}(\phi)$.
 Moreover, since $\mathrm{deg}(M_1)\geq 0$, $\mathrm{deg}(M_2)\geq 0$, 
 $\Phi_i(E_{\geq 0}\otimes_X M_i) \subseteq E_{\geq 0}$. 
It follows that the data 
\[
\CE_{\geq 0}=(E_{\geq 0}, V\otimes \CO_X,\Phi_i|_{E_{\geq 0}\otimes_X M_i}, 0, \psi)
\] 
is a nontrivial subobject of $\CE$. 
If $E_{<0}$ is not the zero sheaf, $\CE_{\geq 0}$ is a proper subobject of 
$\CE$. Then $\delta$-semistability condition implies  
$r(\CE_{\geq 0}) < r(\CE)$, hence 
\be\label{eq:posdegslope}
{d(\CE_{\geq 0})+ v(\CE_{\geq 0})\, \delta \over r(\CE_{\geq 0}) }
 \leq {e+v\, \delta \over r}.
\ee
However $e< d(\CE_{\geq 0})$ and $0<r(\CE_{\geq 0})<r$ under the current assumptions. 
Since also $v(\CE_{\geq 0})=v$ and $\delta, d(\CE_{\geq 0})>0$, inequality 
\eqref{eq:posdegslope} leads to a contradiction. 
Therefore $E_{<0}=0$ and $\phi$ must be identically zero. 

\hfill $\Box$

Let $\CC^0_\CX$ be the full abelian subcategory of $\CC_\CX$ 
consisting of ADHM sheaves $\CE$ with $\phi=0$. For any 
$\delta\in \IR$, an object $\CE$ of $\CC^0_\CX$ will be called 
 $\delta$-semistable 
if it is $\delta$-semistable as an object of $\CC_\CX$. 
Note that given 
an object $\CE$ of $\CC^0_\CX$, any subobject $\CE'\subset \CE$ 
must also belong to $\CC^0_\CX$.  In particular all test subobjects in 
definition (\ref{udeltastability}) also belong to $\CC^0_\CX$, and 
one obtains a stability condition on the abelian category 
$\CC^0_\CX$. Then the properties 
of $\delta$-stability and moduli stacks of semistable objects 
in $\CC^0_\CX$ are analogous to those of 
 $\CC_\CX$. In particular for fixed 
$(r,e,v)\in \IZ_{\geq 1}\times \IZ\times \IZ_{\geq 1}$ there are finitely many 
critical stability parameters of type $(r,e,v)$ dividing the real axis into 
stability chambers. The main difference between $\CC^0_\CX$ and 
$\CC_\CX$ is the presence of an empty chamber, as follows. 
\begin{lemm}\label{emptychamber}
For any $(r,e,v)\in \IZ_{\geq 1}\times \IZ\times \IZ_{\geq 1}$ the moduli 
stack of $\delta$-semistable objects of $\CC^0_\CX$ 
of type $(r,e,v)$ is empty if $\delta<0$. 
\end{lemm}

{\it Proof.} Given an ADHM sheaf $\CE=(E,V,\Phi_i,\psi)$ of 
type $(r,e,v)$, it is straightforward to check that for $\delta<0$ the proper nontrivial 
object $(E,0,\Phi_i,0)$ is always destabilizing if $\delta<0$. 

\hfill $\Box$

\begin{lemm}\label{posdeg} 
Let $\CE$ be a $\delta$-semistable object of $\CC^0_\CX$ of type 
$(r,e,v)\in \IZ_{\geq 1}\times \IZ\times \IZ_{\geq 0}$ for some 
$\delta\geq 0$. 
If $e\geq 0$, then $E_{<0}=0$ and $\phi$ is identically zero. 
\end{lemm}

{\it Proof.} For $\delta>0$ and $v>0$, this obviously follows from lemma (\ref{zerophi}). 
If $\delta=0$ or $v=0$ note that $E_{\geq 0}$ cannot be the zero sheaf 
since $e\geq 0$. Then the proof of lemma (\ref{zerophi}) also applies to this 
case as well. 

\hfill $\Box$

\begin{lemm}\label{genzeroHiggs} 
Let $\CE=(E,0,\Phi_i,0,0)$ be a semistable object of $\CC^0_\CX$
of type $(r,e,0)$, $(r,e)\in \IZ_{\geq 1}\times \IZ$. 
If $(d_1,d_2)=(1,1)$, $E$ must be isomorphic to 
$\CO_X(n)^{\oplus r}$ for some $n\in \IZ$, and $\Phi_i=0$ for $i=1,2$.
If $(d_1,d_2)=(0,2)$, $E$ must be isomorphic to $\CO_X(n)^{\oplus r}$ for some 
$n\in \IZ$, and $\Phi_2=0$. 
\end{lemm}

{\it Proof.}  In both cases, let $E\simeq \oplus_{s=1}^r \CO_X(n_s)$ 
for some $n_s\in \IZ$ so that $n_1\leq n_2\leq \cdots \leq n_r$. 
Since $d_1,d_2\geq 0$, any subsheaf of the form 
\[
\oplus_{s=s_0}^r \CO_X(n_s) 
\]
for some $1\leq s_0\leq r$ must be $\Phi_i$-invariant, 
$i=1,2$.
Therefore the semistability condition implies 
\[
{n_{s_0} +\cdots +n_r\over r-s_0+1} \leq {n_1+\ldots + n_r\over r}
\]
for any $1\leq s_0\leq r$. Then it is straightforward to check that 
$n_1=\cdots = n_r =n$. The rest is obvious. 

\hfill $\Box$.  

\begin{coro}\label{genzeroHiggsinv} 
Under the same conditions as in lemma (\ref{genzeroHiggs}), 
\be\label{eq:Higgsinv}
H(r,e) = \left\{\begin{array}{ll}
{(-1)^{d_1-1}\over r^2} & \mathrm{if}\ e=rn, \ n\in \IZ\\
& \\
0 & {\mathrm{otherwise}}. \\ \end{array}\right. 
\ee
\end{coro}

{\it Proof.} If $(d_1,d_2)=(1,1)$, lemma (\ref{genzeroHiggs}) 
implies that the moduli stack 
$\mfm^{ss}(\CX,r,e,0)$ is isomorphic to the quotient stack 
$[*/GL(r)]$ if $e=rn$ for some $n\in \IZ$, and empty otherwise. 
Alternatively, if $e=rn$, the moduli stack $\mfm^{ss}(\CX,r,e,0)$
can be identified with the moduli stack of trivially semistable 
representations of dimension $r$ of a quiver consisting of only one 
vertex and no arrows. Recall that the trivial semistability condition 
for quiver representations is King stability with all stability 
parameters associated to the vertices set to zero \cite[Ex. 7.3]{genDTI}. 

If $(d_1,d_2)=(0,2)$, lemma (\ref{genzeroHiggs}) 
implies that the moduli stack 
$\mfm^{ss}(\CX,r,rn,0)$, $n\in \IZ$, 
is isomorphic to the moduli stack 
of trivially semistable representations of dimension $r$ of a quiver 
consisting of one vertex and one arrow joining the unique 
vertex with itself. If $e$ is not a multiple of $r$, the moduli stack $\mfm^{ss}(\CX,r,e,0)$
is empty. 

Then 
corollary (\ref{genzeroHiggsinv}) follows 
by a computation very similar to \cite[Sect. 7.5.1]{genDTI}. 

\hfill $\Box$

\begin{rema}\label{reqzerorem}
The same arguments as in the proof of corollary 
(\ref{genzeroHiggsinv}) imply that for any $\delta>0$, 
\be\label{eq:purevinv} 
A_\delta(0,0,1)=1\qquad A_\delta(0,0,2)={1\over 4}.
\ee
\end{rema} 

Extension groups in $\CC^0_\CX$ can be determined by analogy with 
those of $\CC_\CX$. Given two locally free objects $\CE'',\CE'$ of 
$\CC^0_\CX$,  let $\wCC(\CE'',\CE')$ be 
the three term complex of locally free $\CO_X$-modules 
\be\label{eq:zerocomplex}
\begin{aligned} 
0 \to \begin{array}{c} \lochom_{X}(E'',E') \\ \oplus \\ 
\lochom_{X}(V''\otimes \CO_X,V'\otimes \CO_X) \\ \end{array}
& {\buildrel d_1\over \longto} 
\begin{array}{c}  \lochom_{X}(E''\otimes_{X}M_1 ,E') \\ \oplus \\
\lochom_{X}(E''\otimes_{X} M_2, E') \\ \oplus \\ 
\lochom_{X}(V''\otimes \CO_X,E') \\ \end{array} 
 {\buildrel d_2\over \longto} 
\lochom_{X}(E''\otimes_{X}M,E') \to 0 \\
\end{aligned}
\ee
where 
\[
\begin{aligned}
d_1(\alpha,f) = 
(& -\alpha \circ \Phi_{1}'' +\Phi_{1}'\circ (\alpha\otimes 1_{M_1}), 
-\alpha \circ \Phi_{2}''+\Phi_{2}'\circ (\alpha\otimes 1_{M_2}),\\ 
&
-\alpha\circ \psi'' +\psi' \circ f)\\
\end{aligned}  
\]
for any local sections $(\alpha,f)$ of 
the first term 
and 
\[
\begin{aligned}
d_2(\beta_1,\beta_2,\gamma) = &
\beta_1 \circ (\Phi''_2\otimes 1_{M_1}) -
\Phi_{2}'\circ (\beta_1\otimes 1_{M_2}) 
- \beta_2\circ (\Phi''_{1}\otimes 1_{M_2})\\
 & + \Phi_{1}'\circ (\beta_2\otimes 1_{M_1}) \\
\end{aligned}
\]
for any local sections $(\beta_1,\beta_2,\gamma)$ 
of the middle term. The degrees of the three terms in 
\eqref{eq:hypercohA} are $0,1,2$ respectively. 
By analogy with lemma(\ref{extlemma}), the following holds. 
\begin{lemm}\label{zeroexts}
Under the current assumptions, 
$\mathrm{Ext}^k_{\CC^0_\CX}(\CE'',\CE') \simeq 
\IH^k(\wCC(\CE'',\CE'))$ for $k=0,1$. 
\end{lemm}

\begin{lemm}\label{zeroeuler} 
Let $\CE', \CE''$ be two nontrivial locally free objects of $\CC_\CX^0$
of types $(r',e',v'), (r'',e'',v'')\in \IZ_{\geq 1}\times \IZ \times \IZ_{\geq 0}$. 
Suppose that $E'_{<0}=0$, $E''_{<0}=0$ 
for both underlying locally free sheaves $E',E''$. Then 
\be\label{eq:zeroeuler} 
\bal
& \mathrm{dim}(\mathrm{Ext}^0_{\CC_\CX}(\CE'',\CE')) -
\mathrm{dim}(\mathrm{Ext}^1_{\CC_\CX}(\CE'',\CE')) - 
\mathrm{dim}(\mathrm{Ext}^0_{\CC_\CX}(\CE',\CE'')) \\  & +
\mathrm{dim}(\mathrm{Ext}^1_{\CC_\CX}(\CE',\CE'')) =
v'(e''+r'')-v''(e'+r'). \\
\eal
\ee
\end{lemm}

{\it Proof.} 
Note that the complex \eqref{eq:zerocomplex} can be written as the cone of a morphism 
of locally free complexes on $X$
\[
\varrho : \CH[-1] \to \CV
\]
where $\CH$ is the complex obtained from $\wCC(\CE'',\CE')$ by 
omitting all direct summands depending on $V',V''$ (as well as making some 
obvious changes of signs),  and $\CV$ is 
the two term complex 
\[
\bal 
\lochom_{X}(V''\otimes \CO_X,V'\otimes \CO_X) & {\buildrel {}\over\longto}
\lochom_{X}(V''\otimes \CO_X,E')\\
f & \longto \psi'\circ f\\
\eal 
\]
with degrees $0,1$. The morphism $\varrho$ is determined by the 
map 
\[
\bal
\lochom_{X}(E'',E') & {\buildrel {}\over\longto}
\lochom_{X}(V''\otimes \CO_X,E')\\
\alpha &\longto -\alpha \circ \psi''\\
\eal 
\]
Therefore there is a long exact sequence of hypercohomology groups 
\be\label{eq:zeroexseq} 
\xymatrix{ 
0\ar[r] & \IH^0(\CV) \ar[r] & \mathrm{Ext}^0_{\CC^0_{\CX}}(\CE'',\CE')\ar[r] & 
\IH^0(\CH(\CE'',\CE')) \\
\ar[r] & \IH^1(\CV) \ar[r] &\mathrm{Ext}^1_{\CC^0_{\CX}}(\CE'',\CE')\ar[r] & 
\IH^1(\CH(\CE'',\CE')) \\
\ar[r] & \IH^2(\CV) \ar[r] & \cdots & \\}
\ee
Since $E'_{<0}=0$ and $X$ is rational, $\IH^2(\CV) =0$. Obviously, there is 
a similar exact sequence with $\CE',\CE''$ interchanged. Then equation 
\eqref{eq:zeroeuler}  easily follows observing that 
\[
\IH^k(\CH(\CE'',\CE')) \simeq \IH^{3-k}(\IH(\CE', \CE''))^\vee 
\]
for all $0\leq k \leq 3$. 

\hfill $\Box$

{\it Proof of Corollary (\ref{closedform}).}
Let $\gamma=(r,e,v)\in \IZ_{\geq 1}\times \IZ_{\geq 0}\times \IZ_{\geq 0}$ 
be an arbitrary numerical type, and $\delta \in \IR_{\geq 0}$. 
Given any decomposition $\gamma=\gamma_1+\cdots + \gamma_l$, $l\geq 1$
so that 
\[
{e_1+v_1\delta\over r_1}=\cdots = {e_l+v_l\delta\over r_l}=
{e+v\delta\over r}
\]
it is obvious that if $v_i=0$ for some $1\leq i\leq l$ then $e_i\geq 0$. 
Moreover, if $\delta=0$, then $e_i\geq 0$ for all $1\leq i\leq l$. 
In particular this holds for all terms in the right hand side of the defining 
equation of log stack functions \eqref{eq:logstfctB}. 
It also holds for all possible numerical types of Harder-Narasimhan 
filtrations associated to a critical stability parameter $\delta_c\geq 0$ 
as in lemma (\ref{HNfiltrations}). Note that if $\delta_c=0$, the last quotient 
$\CF_h$ in the Harder-Narasimhan filtration with respect to $\delta_+$-stability, 
respectively the first quotient $\CF_1$  in the Harder-Narasimhan filtration with respect 
to $\delta_-$-stability is allowed to be isomorphic to the object $O_v=(0,\IC^v,0,0,0)$, 
$v\geq 1$. 
In conclusion, the definition of generalized Donaldson-Thomas invariants, and 
derivation of wallcrossing formulas carry over to the present set-up for semistable 
objects of positive degree and stability parameters $\delta\geq 0$. 
In this case the resulting invariants will be denoted by $A^0_\delta(\gamma)$, 
or $A^0_\delta(\alpha, v)$ by analogy with section (\ref{wallformula}).
Lemmas (\ref{zerophi}) and (\ref{posdeg}) imply that the invariants $A^0_\delta(\alpha,2)$ 
satisfy the wallcrossing formula \eqref{eq:wallformulaG} at a positive critical 
stability parameter $\delta_c$ of type $(\alpha,2)$. 
If $\delta_c=0$, a modification of formula \eqref{eq:wallformulaG} 
is required, reflecting the presence of objects isomorphic to $\CO_v$, $v=1,2$ 
in the Harder-Narasimhan filtrations. Basically one has to set $\delta_c=0$ 
in conditions \eqref{eq:slopesA}-\eqref{eq:slopesC}, and 
allow elements $(\alpha_i)_{1\leq i\leq l}$ so that $\alpha_i$, $1\leq i\leq l-1$ 
satisfy conditions \eqref{eq:slopesA}-\eqref{eq:slopesC}, and $\alpha_l=(0,0)$. 
This will result in extra terms in the right hand side of \eqref{eq:wallformulaG} 
which can be easily written down using \eqref{eq:purevinv}.
Since this is an easy exercise, explicit formulas will be omitted (see 
\cite[Thm. 1.$ii$.]{chamberII} for the $v=1$ case). 
Finally, note that one can also check compatibility with the Kontsevich-Soibelman 
formula at $\delta_c=0$ repeating the calculations in section (\ref{KSsect}). 

Then the proof of corollary (\ref{closedform}) will be based on the 
KS wallcrossing formula relating $\delta$-invariants for $\delta<0$ 
to $\delta$-invariants with $\delta>>0$. 
Let $(r,e)\in \IZ_{\geq 1}\times \IZ_{\geq 0}$
and let $\delta_+\in \IR_{>0}\setminus \IQ$ an irrational 
stability parameter so that $\delta_+$ is asymptotic 
of type $(r',e')$ for all $1\leq r'\leq r$, $0\leq e'\leq e$, $1\leq v\leq 2$. 
Moreover, assume that $re<\delta_+$. 
Then the KS formula reads 
\be\label{eq:emptytoinftyA} 
\bal 
\prod_{(r,n, v) \in \IZ_{\geq 1}\times\IZ_{\geq 0}\times \{0,1,2\} \cup \{0,0,1\} }
U_{\lambda(r,n,v)}^{{\overline {A^0}}_-(r,e,v)} = 
\prod_{(r,n, v) \in \IZ_{\geq 1}\times \IZ_{\geq 0}\times \{0,1,2\} \cup \{0,0,1\}  }
U_{\lambda(r,n,v)}^{{\overline {A^0}}_+(r,n,v)} 
\eal 
\ee
where in each term the factors are ordered in increasing order of $\delta_\pm$-slopes 
from left to right. 
The alternative notation introduced in section (\ref{KSsect}) will be used in the 
following. 
Then 
corollary (\ref{genzeroHiggsinv}) and equation \eqref{eq:purevinv} imply that 
the left hand side of \eqref{eq:emptytoinftyA} reads 
\be\label{eq:emptytoinftyB} 
\mathrm{exp}({ {\sf f}}_{00} + \frac{1}{4}{{\sf g}}_{00} )
\prod_{n=0}^{\infty}U_{{{\sf e}}_{1n}}, 
\ee
where 
\[
U_{{{\sf e}}_{rn}} = \mathrm{exp}\bigg((-1)^{d_1-1}\sum_{k=1}^{\infty} \frac{
{ {\sf e}}_{kr,kn}}{k^2}\bigg).
\]
Moreover, given the above choice of $\delta_+$, 
\[
e< {\delta_+ \over r} < \cdots {e+\delta_+\over r} < {\delta_+\over r-1} < \cdots 
< {e+\delta_+\over r-1} < \cdots < \delta_+ + e < {2\delta_+ \over r} < \cdots < 
2\delta_+ + e. 
\]
Therefore, 
in the right hand side of equation \eqref{eq:emptytoinftyA}, 
the factors of the form $U_{{\lambda}(r',e',v)}^{{\overline A}_+(r',e',v)}$, 
with $v\in \{0,1,2\}$, 
and $1\leq r'\leq r$, $1\leq e'\leq e$ occur in the following 
order 
\be\label{eq:emptytoinftyC} 
\bal 
& \prod_{n=0}^{e}U_{{{\sf e}}_{1n}} 
\prod_{n=0}^e U_{{{\sf f}}_{r,n}}^{{\overline A}_+(r,n,1)}
\prod_{n=0}^e U_{{{\sf f}}_{r-1,n}}^{{\overline A}_+(r-1,n,1)}\cdots 
\prod_{n=0}^e U_{{{\sf f}}_{1,n}}^{{\overline A}_+(1,n,1)} U_{{{\sf f}}_{0,0}}^{{\overline A}_+(0,0,1)}\\
& 
\prod_{n=0}^e U_{{{\sf g}}_{r,n}}^{{\overline A}_+(r,n,2)}\cdots 
\prod_{n=0}^e U_{{{\sf g}}_{r-1,n}}^{{\overline A}_+(r-1,n,2)}\cdots 
\prod_{n=0}^e U_{{{\sf g}}_{1,n}}^{{\overline A}_+(1,n,2)},\\
\eal 
\ee
where 
\[
U_{{{\sf f}}_{rn}} = \mathrm{exp}({ {\sf f}}_{rn} + 
\frac{1}{4} {{\sf g}}_{2r,2n}),\qquad 
U_{{ {\sf g}}_{rn}} = \mathrm{exp}({ {\sf g}}_{rn}).
\]
In addition, the right hand side of  \eqref{eq:emptytoinftyA} contains of course 
extra factors of the form $U_{{\lambda}(r',e',v)}^{{\overline A}_+(r',e',v)}$,
with $v\in \{0,1,2\}$, 
and either $ r'> r$ or  $e'> e$. Some of these extra factors may in fact occur
between the factors listed in \eqref{eq:emptytoinftyC}. 
However,  they can be ignored for the purpose 
of this computation since commutators involving such factors are again 
expressed in terms of generators ${\lambda}(r',e',v)$ with 
either  $ r'> r$ or  $e'> e$. Therefore, using the BCH formula, 
\eqref{eq:emptytoinftyA} 
yields 
\be\label{eq:KSconifold2}
\begin{aligned}
& (\prod_{n=0}^{e}U_{{{\sf e}}_{1n}})^{-1} \,
\mathrm{exp}({{\sf f}}_{00} + \frac{1}{4}{{\sf g}}_{00} )\,
\prod_{n=0}^{\infty}U_{{{\sf e}}_{1n}}  = \\
& \mathrm{exp} \Big(\mathop{\sum_{1\leq s\leq r,\ 0\leq n \leq e}}_{}
 {A}_{+}(s,n,1)\, { {\sf f}}_{sn}+ 
\mathop{\sum_{ 1\leq s\leq r,\ 0\leq n\leq e}}_{}{ A}_{+}(s,n,2)\, 
{ {\sf g}}_{sn} +\\
& \sum_{\substack{r_1>r_2\geq 1, \ r_1+r_2\leq r,\ 
n_1,\ n_2\geq 0, n_1+n_2\leq e \\ \mathrm{or}  \ 
1\leq r_1=r_2\leq r/2, \ 0\leq n_1< n_2, \ n_1+n_2\leq e \\ \mathrm{or}  \ 1 \leq r_1 \leq r, \ 0 \leq n_1 \leq
e, \ r_2=n_2=0}}  \frac{1}{2} (n_1-n_2+r_1-r_2) (-1)^{(n_1-n_2+r_1-r_2)} \\
& {A}_{+}
(r_1,n_1,1){ A}_{+}(r_2,n_2,1)\,
{ {\sf g}}_{r_1+r_2,n_1+n_2}+\cdots\Big)\\
\end{aligned}
\ee
where $\cdots$ are terms involving generators ${\lambda}(r',e',v)$ with 
either  $ r'> r$ or  $e'> e$. 
For fixed $e\geq 1$, let $\mathcal{H}_e$ be defined by
\be
\mathrm{exp}(\mathcal{H}_e)\equiv \prod_{n=0}^{e}U_{{\sf e}_{1,n}} = 
\mathrm{exp} \bigg( (-1)^{d_1-1}\sum_{ 0\leq n \leq e ,\ k\geq 1}   
\frac{{\sf e}_{k,kn}}{k^2} \bigg).
\ee
Using the BCH formula, the left hand side of equation \eqref{eq:KSconifold2} 
becomes 
\be\label{eq:KSconifold3} 
\bal 
\mathrm{exp}\bigg({\sf f}_{00} +{1\over 4}{\sf g}_{00} 
+ \sum_{j=1}^\infty {1\over j!} [\underbrace{-\mathcal{H}_e, \cdots [ -\mathcal{H}_e}_{\text{j times}},{\sf f}_{00} +{1\over 4}{\sf g}_{00} ]\cdots ]\bigg) 
\eal 
\ee
modulo terms involving generators ${\lambda}(r',e',v)$ with 
either  $ r'> r$ or  $e'> e$. 

Next, the  Lie algebra commutators 
\[
\bal 
& [{\sf e}_{r_1,n_1}, {\sf f}_{r_2,n_2}] = (-1)^{n_1+r_1}(n_1+r_1)\,
{\sf f}_{r_1+r_2,n_1+n_2}\\
& [{\sf e}_{r_1,n_1}, {\sf g}_{r_2,n_2}] = 2(n_1+r_1)\, 
{\sf g}_{r_1+r_2,n_1+n_2},\\
\eal 
\]
yield
\[
\bal 
& [ \underbrace{-\mathcal{H}_e, \cdots [ -\mathcal{H}_e}_{\text{j times}},{\sf f}_{00}
 ]\cdots ] =
\mathop{\sum_{n_1,\ldots, n_j = 0}^{e}}_{} 
\mathop{\sum_{k_1,\ldots, k_j \geq 1}}_{} 
(-1)^{j(d_1-1)}
\prod_{i=1}^j {n_i+1\over k_i} (-1)^{(n_i+1)k_i-1} \, {\sf f}_{k_1+\cdots + k_j, 
k_1n_1+\cdots + k_jn_j}\\
\eal 
\]
and 
\[
\bal 
& [ \underbrace{-\mathcal{H}_e, \cdots [ -\mathcal{H}_e}_{\text{j times}},{\sf g}_{00} 
]\cdots ] =
\mathop{\sum_{n_1,\ldots, n_j = 0}^{e}}_{} 
\mathop{\sum_{k_1,\ldots, k_j \geq 1}}_{} 
(-1)^{j(d_1-1)}\prod_{i=1}^j (-2){n_i+1\over k_i}  \, {\sf g}_{k_1+\cdots + k_j, 
k_1n_1+\cdots + k_jn_j}\\
\eal 
\]
Therefore, identifying the coefficients of the generators 
${\sf f}_{rn}$ in \eqref{eq:KSconifold2}
it follows that the invariant $A_+(r',e',1)$ with $1\leq r'\leq r$ and 
$0\leq e'\leq e$ equals  the coefficient of the monomial $u^{r'}q^{e'+r'}$ 
in the expression 
\[
\bal 
\sum_{j=0}^\infty {1\over j!} 
\left(
\mathrm{ln}\left(\prod_{n=0}^e (1-u(-q)^{n+1})^{(-1)^{d_1-1}(n+1)}\right)\right)^j =
\prod_{n=1}^{e+1}(1-u(-q)^n)^{(-1)^{d_1-1}n}.
\eal
\]
Similarly, identifying the coefficients of the generators ${\sf g}_{rn}$ in \eqref{eq:KSconifold2}
proves that the invariant $A_+(r',e',2)$ with $1\leq r'\leq r$ and 
$0\leq e'\leq e$ equals  the coefficient of the monomial $u^{r'} q^{e'+r'}$ 
in the expression
\[
\bal
& {1\over 4}\prod_{n=1}^{e+1}(1-uq^n)^{2(-1)^{d_1-1}n}
- 
\sum_{\substack{r_1>r_2\geq 1, \ r_1+r_2\leq r,\ 
n_1,\ n_2\geq 0, n_1+n_2\leq e \\ \mathrm{or}  \ 
1\leq r_1=r_2\leq r/2, \ 0\leq n_1< n_2, \ n_1+n_2\leq e 
\\ \mathrm{or}  \ 1 \leq r_1 \leq r, \ 0 \leq n_1 \leq
e, \ r_2=n_2=0}} \\
& \frac{1}{2}  (n_1+r_1-n_2-r_2) (-1)^{(n_1+r_1-n_2-r_2)} 
A_{+}(r_1,n_1,1)
A_{+}(r_2,n_2,1) q^{r_1+r_2}u^{n_1+n_2}.\\
\end{aligned}
\]
Since this holds for any $(r,e)\in \IZ_{\geq 1}\times \IZ_{\geq 0}$
(with a suitable choice of $\delta_+$), corollary (\ref{closedform}) 
follows. 

\hfill $\Box$

\bibliography{adhmref.bib}
 \bibliographystyle{abbrv}
\end{document}